\documentclass[11pt,twoside]{article}

\usepackage{amssymb}
\usepackage{amsmath}

\allowdisplaybreaks

\def\rr{{\mathbb R}}
\def\rd{{{\rr}^d}}

\def\zz{{\mathbb Z}}
\def\nn{{\mathbb N}}
\def\hh{{\mathbb H}}
\def\bbg{{\mathbb G}}

\def\cc{{\mathbb C}}

\def\cx{{\mathcal X}}

\def\cd{{\mathcal D}}

\def\cl{{\mathcal L}}

\def\cb{{\mathcal B}}

\def\fz{\infty}
\def\az{\alpha}

\def\loc{{\mathop\mathrm{\,loc\,}}}

\def\lz{\lambda}
\def\dz{\delta}

\def\ez{\epsilon}
\def\ezl{{\epsilon_1}}

\def\kz{\kappa}
\def\bz{\beta}

\def\gz{{\gamma}}

\def\wz{\widetilde}

\def\hs{\hspace{0.3cm}}

\def\ls{\lesssim}
\def\gs{\gtrsim}

\def\hl{{H^1_\cl(\rn)}}

\def\lt{{L^2(\cx)}}
\def\lin{{L^\fz(\cx)}}

\def\hl{{\mathop\mathrm {HL}}}

\def\ati{{\mathrm {AOTI}}}
\def\bmo{{\mathop\mathrm{BMO}}}
\def\blo{{\mathop\mathrm{BLO}}}
\def\bmox{{\mathop\mathrm{BMO}(\cx)}}
\def\blox{{\mathop\mathrm{BLO}(\cx)}}
\def\bmoz{{\mathrm{BMO}_\rho(\cx)}}
\def\bloz{{\mathrm{BLO}_\rho(\cx)}}

\def\infi{{\mathop\mathrm{essinf}}}

\def\einf{{\mathop{\mathrm{\,essinf\,}}}}
\def\dsum{\displaystyle\sum}
\def\diam{{\mathop\mathrm{\,diam\,}}}

\def\dint{\displaystyle\int}

\def\dfrac{\displaystyle\frac}
\def\dsup{\displaystyle\sup}

\def\r{\right}
\def\lf{\left}

\newtheorem{thm}{Theorem}[section]
\newtheorem{lem}{Lemma}[section]
\newtheorem{prop}{Proposition}[section]
\newtheorem{rem}{Remark}[section]
\newtheorem{cor}{Corollary}[section]
\newtheorem{defn}{Definition}[section]

\newtheorem{pf}{\it Proof.}[section]

\begin{document}

\arraycolsep=1pt

\title{{\vspace{-5cm}\small\hfill\bf Commun. Pure Appl. Anal., to appear}\\
\vspace{4cm}\Large\bf LOCALIZED BMO AND BLO SPACES ON ${\rm RD}$-SPACES
AND APPLICATIONS TO SCHR\"ODINGER OPERATORS
\footnotetext{\hspace{-0.35cm} 2000 {\it Mathematics Subject
Classification}. {Primary 42B35; Secondary 42B20, 42B25, 42B30.}
\endgraf{\it Key words and phrases.} Space of homogeneous type,
Heisenberg group, connected and simply connected nilpotent Lie group, reverse
H\"older inequality, admissible function, Schr\"odinger operator,
localized $\bmo$ space, localized $\blo$ space, maximal function, $g$-function.
\endgraf
 The first author is supported by the National
Natural Science Foundation (Grant No. 10871025) of China.}}
\author{\sc{Dachun Yang, Dongyong Yang and Yuan Zhou}\\
\footnotesize School of Mathematical Sciences,
 Beijing Normal University, \\
 \footnotesize Laboratory of Mathematics and Complex systems,
Ministry of Education,\\
\footnotesize Beijing 100875, People's Republic of China}
\date{ }
\maketitle

\begin{center}
\begin{minipage}{13.5cm}\small
{\noindent{\sc Abstract.}
An RD-space ${\mathcal X}$ is a space of homogeneous
type in the sense of Coifman and Weiss with the additional property
that a reverse doubling condition holds in ${\mathcal X}$.
Let $\rho$ be an admissible function on RD-space ${\mathcal X}$.
The authors first introduce the localized spaces
$\mathrm{\,BMO}_\rho({\mathcal X})$ and $\mathrm{\,BLO}_\rho({\mathcal X})$
and establish their basic properties, including the John-Nirenberg
inequality for $\mathrm{\,BMO}_\rho({\mathcal X})$,
several equivalent characterizations for $\mathrm{\,BLO}_\rho({\mathcal X})$,
and some relations between these spaces.
Then the authors obtain the boundedness on these localized spaces
of several operators including the natural maximal operator,
the Hardy-Littlewood maximal operator,
the radial maximal functions and their localized
versions associated to $\rho$, and the Littlewood-Paley $g$-function
associated to $\rho$, where the Littlewood-Paley $g$-function
and some of the radial maximal functions are
defined via kernels which are modeled on the semigroup
generated by the Schr\"odinger operator.
These results apply in a wide range of settings, for instance,
to the Schr\"odinger operator or
the degenerate Schr\"odinger operator on ${{\mathbb R}}^d$,
or the sub-Laplace Schr\"odinger operator on Heisenberg
groups or connected and simply connected nilpotent Lie groups.}
\end{minipage}
\end{center}

\vspace{-1cm}
\section*{}\label{s1}
\noindent {\bf 1. Introduction.}
 Since the space, $\bmo(\rd)$, of functions with
bounded mean oscillation on $\rd$ was introduced by
John and Nirenberg \cite{jn}, it then plays an important
role in harmonic analysis and partial differential equations.
For example, it is well known that $\bmo(\rd)$ is the dual space
of the Hardy space $H^1(\rd)$ (see, for example, \cite{s93,g}),
and also a good substitute of $L^\fz(\rd)$.
Recall that the Riesz transforms $\nabla(-\Delta)^{-1/2}$ are bounded on
$\bmo(\rd)$ but not on $L^\fz(\rd)$ (see again, for example, \cite{s93,g}),
where $\Delta\equiv\sum_{j=1}^d\frac{\partial^2}{\partial x_j^2}$
is the Laplacian and $\nabla$ is the gradient operator.
However, the space $\bmo(\rd)$ is essentially
related to the Laplacian $\Delta$.

Let $\cl\equiv -\Delta+V$ be the Schr\"odinger operator on $\rd$,
where the potential $V$ is a nonnegative locally integrable function.
Recently, there is an increasing interest on the study of these
operators. In particular, Fefferman \cite{f83},
Shen \cite{s95} and Zhong \cite{z99} established
some basic results, including some estimates of the fundamental solutions and the
boundedness on Lebesgue spaces of Riesz transforms, for $\cl$ on $\rd$
with $d\ge3$ and the nonnegative potential $V$ satisfying
the reverse H\"older inequality. Especially, the works of Shen \cite{s95}
lay the foundation for developing harmonic analysis related to
$\cl$ on $\rd$. Li \cite{l99} extended part of these
results in \cite{s95} to the sub-Laplace Schr\"odinger operator
on connected and simply connected nilpotent Lie groups.
On the other hand, denote by $\cb_q(\rd)$ the class of functions
satisfying the reverse H\"older inequality of order $q$.
For $V\in\cb_{d/2}(\rd)$ with $d\ge3$,
Dziuba\'nski et al \cite{dgmtz05}
introduced the BMO-type space $\bmo_\cl(\rd)$
associated to the auxiliary function $\rho$ determined by the potential $V$
(see, for example, \eqref{2.4} below)
and established the duality between $H^1_\cl(\rd)$
and $\bmo_\cl(\rd)$,
as well as a characterization of $\bmo_\cl(\rd)$ in
terms of the Carleson measure and
the $\bmo_\cl(\rd)$ boundedness of the variants of some classical operators
associated to $\cl$ including semigroup
maximal functions and the Hardy-Littlewood maximal function.
These results were generalized to Heisenberg groups by
Lin and Liu \cite{ll08}. Also, it is now known that $\bmo_\cl(\rd)$
in \cite{dgmtz05} is a special case of BMO-type spaces introduced
by Duong and Yan \cite{dy1,dy2}; see, in particular, \cite[Proposition
6.11]{dy2} and also \cite{yz}.

Recently, a theory of Hardy spaces and their dual spaces
on so-called RD-spaces was established in \cite{hmy1,hmy2,gly1}.
A space of homogenous type $\cx$ in the sense of Coifman and Weiss
(\cite{cw71,cw77}) is called
to be an RD-space if $\cx$ has the additional property that a
reverse doubling condition holds in $\cx$ (see \cite{hmy2}). It is
well known that a connected space of homogeneous type is an
RD-space. Typical examples of RD-spaces include Euclidean spaces,
Euclidean spaces with weighted measures satisfying the
 doubling property, Heisenberg groups, connected and simply
 connected nilpotent Lie groups (\cite{v88,v92}) and the boundary of an unbounded
model polynomial domain in $\cc^N$ (\cite{ns06}), or more
generally, Carnot-Carath\'eodory spaces with doubling measures
(\cite{nsw85, hmy2}). In \cite{yz08}, modeled on the known auxiliary function
determined by $V$, a notion of admissible
functions $\rho$ was introduced and a theory of the localized
Hardy space $H^1_\rho(\cx)$ associated with a given admissible function
$\rho$ was developed. In particular, the space $H^1_\rho(\cx)$ was
characterized via several
maximal functions modeled on the semigroup maximal operators
generated by Schr\"odinger operators, including the localized
radial maximal function $S^+_\rho$.

One of the main purposes of this paper is to
investigate behaviors of these maximal operators aforementioned
on localized  BMO spaces.
Precisely, let $\rho$ be an admissible function on RD-space $\cx$.
We first introduce the localized BMO space
$\bmo_\rho(\cx)$ and localized BLO space $\blo_\rho(\cx)$,
and establish their basic properties,
including the John-Nirenberg inequality for $\bmo_\rho(\cx)$,
several equivalent characterizations for $\blo_\rho(\cx)$,
and some relations between these spaces.
Then we obtain the boundedness on these localized spaces
of several operators including the natural maximal operator,
the Hardy-Littlewood maximal operator,
the radial maximal functions and their localized
versions associated to $\rho$, and the Littlewood-Paley $g$-function
associated to $\rho$, where the Littlewood-Paley $g$-function
and some of the radial maximal functions are
defined via kernels which are modeled on the semigroup
generated by the Schr\"odinger operator.
These results apply in a wide range of settings.
Moreover, even when these results are applied, respectively, to
the Schr\"odinger operator or
the degenerate Schr\"odinger operator on $\rd$,
or the sub-Laplace Schr\"odinger operator on Heisenberg
groups or connected and simply connected nilpotent Lie groups,
we also obtain some new results.

To be precise, this paper is organized as follows.

In Section 2, we first recall some notation and
notions from \cite{hmy2,yz08},
including the approximation of the identity, the admissible function $\rho$,
the radial maximal function $S^+(f)$
and the localized radial maximal function $S^+_\rho(f)$,
where $S^+(f)$ and $S^+_\rho(f)$ are defined
via a given approximation of the identity.

In Section 3, letting $\rho$ be an admissible function on $\cx$,
we first introduce the localized BMO space
$\bmo_\rho(\cx)$ and localized BLO space $\blo_\rho(\cx)$;
see Definitions \ref{d3.1} and \ref{d3.2} below.
We also recall the notions of their global versions in
Definitions \ref{d3.1} and \ref{d3.2} below.
Then we establish some useful properties concerning these spaces, including the
John-Nirenberg inequality for $\bmo_\rho(\cx)$ (see Theorem \ref{t3.1} below),
several characterizations and inclusion relations of these
spaces (see Lemma \ref{l3.1},
Remarks \ref{r3.1} and \ref{r3.2}, and Corollary \ref{c3.1} below).
Then we prove that the function in $\blo_\rho(\cx)$
has lower bound in Theorem \ref{t3.2},
and establish several equivalent characterizations
of $\blo_\rho(\cx)$ in Theorems \ref{t3.2} and \ref{t3.3}, Remark \ref{r3.3},
and Corollaries \ref{c3.2} and \ref{c3.3} below.

In Section 4, we establish the boundedness of
the natural maximal function, the Hardy-Littlewood
maximal function and their localized versions from
$\bmoz$ to $\bloz$, and as an application, we obtain
several equivalent characterizations for $\blo_\rho(\cx)$
via the localized natural maximal function; see Theorems \ref{t4.1} and
\ref{t4.2}, Lemma \ref{l4.1} and Corollary \ref{c4.1} below.
We point out that Corollary \ref{c4.1} improves the
results of \cite{dgmtz05} and \cite{ll08}
even for the Schr\"odinger operators on $\rd$ or
Heisenberg groups with the potentials
satisfying certain reverse H\"older inequality; see Remark \ref{r4.1} below.

In Section 5, we establish the boundedness of some
maximal operators from $\bmoz$ to $\bloz$.
To be precise, the boundedness
of the radial maximal functions $S^+(f)$, $S^+_\rho(f)$
and certain maximal operator $T^+$ from $\bmoz$ to $\bloz$ are
presented in Section 5.1; see Theorem \ref{t5.1},
Corollaries \ref{c5.1} and \ref{c5.2} below.
These operators were used, respectively in \cite{gly1} and \cite{yz08},
to characterize the corresponding  Hardy spaces $H^1(\cx)$ and $H^1_\rho(\cx)$.
Section 5.2 is devoted to the boundedness of
$P^+$ from $\bmoz$ to $\bloz$; see
Theorem \ref{t5.2} below. Here, $T^+$ and $P^+$
are defined via kernels which are modeled on the semigroup
generated by the Schr\"odinger operator, and were used in
\cite{yz08} to characterize the corresponding Hardy space
$H^1_\rho(\cx)$.

In Section 6,
we obtain the boundedness on $\bmo_\rho(\cx)$ of the Littlewood-Paley
$g$-function which is also defined via kernels modeled on the semigroup
generated by the Schr\"odinger operator.
Assuming that $g$-function is bounded on $L^2(\cx)$, we prove that
if $f\in\bmoz$, then $[g(f)]^2\in\bloz$
with norm no more than
$C\|f\|_\bmoz^2$, where $C$ is a
positive constant independent of $f$;
see Theorem \ref{t6.1} below. As a corollary,
we obtain the boundedness of the Littlewood-Paley $g$-function
from $\bmoz$ to $\bloz$;
see Corollary \ref{c6.1} below.

In Section 7, we apply results obtained in Sections 5
and 6, respectively, to the Schr\"odinger operator or the
degenerate Schr\"odinger operator on $\rd$, the
sub-Laplace Schr\"odinger operator on Heisenberg groups or on
connected and simply connected nilpotent Lie groups.
The nonnegative potentials of these Schr\"odinger operators are assumed to
satisfy the reverse H\"older inequality. See Propositions \ref{p7.2},
\ref{p7.3}, \ref{p7.4}
and \ref{p7.5} below.
Even for these special cases, our results
further improve and generalize
the corresponding results in \cite{dgmtz05,ll08}.

We now make some conventions. Throughout this paper, we always use
$C$ or $A$ to denote a positive constant that is independent of the main
parameters involved but whose value may differ from line to line.
Constants with subscripts, such as $C_1$ or $A_1$, do not change in different
occurrences. If $f\le Cg$, we then write $f\ls g$ or $g\gs f$; and
if $f \ls g\ls f$, we then write $f\sim g.$ For any given ``normed"
spaces $\mathcal A$  and $\mathcal B$,
the symbol ${\mathcal A}\subset {\mathcal B}$ means that for all $f\in \mathcal A$,
then $f\in\mathcal B$ and $\|f\|_{\mathcal B}\ls \|f\|_{\mathcal A}$.
We always use $B$ to denote a ball of $\cx$, and for any ball $B\subset
\cx$, we denote by $x_B$ the center of $B$, $r_B$ the radius of $B$,
and $B^\complement\equiv\cx\setminus B$. Moreover, for any ball
$B\subset \cx$ and $\lz>0$, we denote by $\lz B$ the ball centered
at $x_B$ and having radius $\lz r_B$.
 Also, $\chi_E$ denotes the characteristic
 function of any set $E\subset\cx$.
For all $f\in L^1_\loc(\cx)$ and balls $B$, we always set
$f_{B}\equiv\frac1{\mu(B)}\int_Bf(y)\,d\mu(y)$.

\vspace{-1cm}
\section*{}
\setcounter{section}{2}
\setcounter{thm}{0}
\setcounter{defn}{0}
\setcounter{lem}{0}
\setcounter{rem}{0}
\setcounter{prop}{0}
\setcounter{cor}{0}
\noindent 2. {\bf Preliminaries}. We first recall the notions of
spaces of homogeneous type in the
sense of Coifman and Weiss \cite{cw71,cw77} and RD-spaces in
\cite{hmy2}.

\begin{defn}\hspace{-0.2cm}{\bf.}\hspace{0.2cm}\label{d2.1}\rm
Let $(\cx,\, d)$ be a metric space with a regular Borel measure
$\mu$ such that all balls defined by $d$ have finite and positive
measure. For any $x\in \cx $ and $r\in(0, \fz)$, set the ball $B(x,r)\equiv\{y\in
\cx :\ d(x,y)<r\}.$

(i) The triple $(\cx,\,d,\,\mu)$ is called a
    space of homogeneous type if there exists
    a constant $A_1\in[1, \fz)$ such that for all $x\in \cx $ and $r\in(0, \fz)$,
    \begin{equation}\label{2.1}
    \mu(B(x, 2r))\le A_1\mu(B(x,r))\ (\mathrm{\it doubling\ property}).
    \end{equation}

(ii) Let $\kz\in(0, n]$. The triple $(\cx,\,d,\,\mu)$ is called a
    $(\kz,\,n)$-space if there
    exist constants $A_2\in(0, 1]$ and $A_3\in[1, \fz)$ such that for all
    $x\in \cx$, $r\in(0, \diam(\cx)/2]$ and $\lz\in[1, \diam(\cx)/(2r)]$,
     \begin{equation}\label{2.2}
     A_2\lz^\kz\mu(B(x,r))\le\mu(B(x,\lz r))\le A_3\lz^ n\mu(B(x,r)),
        \end{equation}
    where $\diam(\cx)\equiv\sup_{x,\,y\in\cx}d(x,y)$.

    A space of homogeneous type is called an RD-space,
    if it is a $(\kz,\,n)$-space for some $\kz\in(0, n]$, i.\,e., if some
    ``reverse" doubling condition holds.
\end{defn}

Obviously, a $(\kappa,\, n)$-space is a space of homogeneous type
with $A_1=A_32^n$. Conversely, a space of homogeneous type satisfies
the second inequality of \eqref{2.2} with $A_3=A_1$ and $n=\log_2A_1$.
Moreover, it was proved in  \cite[Remark 1]{hmy2} that $\cx$
is an RD-space if and only if
$\cx$ is a space of homogeneous type with the additional property
that there exists a constant $a_0\in(1, \fz)$ such that for all $x\in\cx$ and
$r\in(0, \diam(\cx)/a_0)$, $B(x,\,a_0r)\setminus B(x,\,r)\neq\emptyset$.

In what follows, we always set
$V_r(x)\equiv\mu(B(x,\,r))$ and $V(x,\,y)\equiv\mu(B(x,\,d(x,\,y)))$
for all $x,\,y\in\cx$ and $r\in(0,\,\fz)$.

\begin{defn}\label{d2.2}\hspace{-0.2cm}{\bf.}\hspace{0.2cm}\rm (\cite{yz08})
A positive function $\rho$ on $\cx$ is said to be admissible if there
exist positive constants $C_0$ and $k_0$ such that for all
$x,\,y\in\cx$,
\begin{equation}\label{2.3}
\frac1{\rho(x)} \le
C_0\frac1{\rho(y)}\lf(1+\frac{d(x,\,y)}{\rho(y)}\r)^{k_0}.
\end{equation}
\end{defn}

We remark that the function $\rho$ in Definition \ref{d2.2} does
exist. Obviously, if $\rho$ is a constant function, then $\rho$ is
admissible. Moreover, let $x_0\in\cx$ being fixed.
The function $\rho(y)\equiv (1+d(x_0,\,y))^s$ for all
$y\in\cx$ with $s\in(-\fz, 1)$ also satisfies Definition \ref{d2.2} with
$k_0=s/(1-s)$ when $s\in[0, 1)$ and $k_0=-s$ when $s\in(-\fz, 0)$. Another
non-trivial class of admissible functions is given by the well-known
reverse H\"older class $\mathcal\cb_q(\cx, d, \mu)$, which
is always written as $\mathcal\cb_q(\cx)$.
Recall
that a nonnegative potential $V$ is said to be in $\cb_q(\cx)$
with $q\in(1,\,\fz]$ if there exists a positive
constant $C$ such that for all balls $B$ of $\cx$,
$$\lf(\frac1{\mu(B)}\dint_B[V(y)]^q\,d\mu(y)\r)^{1/q}
\le \frac{C}{\mu(B)}\dint_BV(y)\,d\mu(y)$$
with the usual modification made when $q=\fz$. It is known that if
$V\in \cb_q(\cx)$ for certain $q\in(1,\,\fz]$, then $V$
is an $A_\fz(\cx)$ weight in the sense of Muckenhoupt, and also
$V\in \cb_{q+\ez}(\cx)$ for some $\ez\in(0, \fz)$; see, for
example, \cite{s93} and \cite{st89}.
Thus $\cb_q(\cx)=\cup_{q_1>q}\cb_{q_1}(\cx)$. For all $V\in
\cb_q(\cx)$ with certain $q\in(1,\,\fz]$ and all
$x\in\cx$, set
\begin{equation}\label{2.4}
\rho(x)\equiv[m(x, V)]^{-1}\equiv\sup\lf\{r>0:\hs \frac{r^2}
{\mu(B(x,\,r))}\dint_{B(x,\,r)}V(y)\,d\mu(y)\le 1\r\};
\end{equation}
see, for example, \cite{s95} and also \cite{yz08}.
It was also proved in \cite{yz08} that
$\rho$ in \eqref{2.4} is an admissible function
if $n\ge1$, $q>\max\{1,\,n/2\}$
and $V\in \cb_q(\cx)$.

The following notion of approximations of the identity
on RD-spaces was first introduced in \cite{hmy2}, whose
existence was given in Theorem 2.1 of \cite{hmy2}.

\begin{defn}\hspace{-0.2cm}{\bf.}\hspace{0.2cm}\label{d2.3}\rm
Let $\ez_1\in(0,\,1]$ and $\ez_2\in(0,\,\fz)$. A sequence
$\{S_k\}_{k\in\zz}$ of bounded linear integral operators on
$L^2(\cx)$ is said to be an approximation of the identity of order
$(\ez_1,\,\ez_2)$ (for short, $(\ez_1,\,\ez_2)$-$\ati$), if there
exists a positive constant $A_4$ such that for all $k\in\zz$ and all
$x,\, x',\, y$ and $y'\in\cx$, $S_k(x,y)$, the integral kernel of
$S_k$ is a measurable function from $\cx\times\cx$ into $\cc$
satisfying
\begin{enumerate}
\vspace{-0.3cm}
\item[(i)] $|S_k(x,\,y)|\le A_4\frac 1{V_{2^{-k}}(x)+V_{2^{-k}}(y)+V(x,\,y)}
\frac{2^{-k\ez_2}}{(2^{-k}+d(x,\,y))^{\ez_2}};$ \vspace{-0.3cm}
\item[(ii)] $|S_k(x,\,y)-S_k(x',\,y)|\le A_4
(\frac{d(x,\,x')}{2^{-k}+d(x,\,y)})^\ezl \frac
1{V_{2^{-k}}(x)+V_{2^{-k}}(y)+V(x,\,y)}
\frac{2^{-k\ez_2}}{(2^{-k}+d(x,\,y))^{\ez_2}}$ \newline for
$d(x,x')\le(2^{-k}+d(x,\,y))/2;$ \vspace{-0.3cm}
\item[(iii)] Property (ii) also holds with $x$ and $y$ interchanged;
\vspace{-0.3cm}
\item[(iv)] $\int_\cx S_k(x,z)\,d\mu(z)=1=
\int_\cx S_k(z,y)\,d\mu(z)$ for all $x, \, y\in\cx$.
\end{enumerate}
\end{defn}

\begin{rem}\hspace{-0.2cm}{\bf.}\hspace{0.2cm}\label{r2.1}\rm
If a sequence $\{\wz S_t\}_{t>0}$ of bounded linear integral
operators on $L^2(\cx)$ satisfies (i) through (iv) of Definition
\ref{d2.3} with  $2^{-k}$ replaced by $t$, then we call $\{\wz
S_t\}_{t>0}$ a continuous approximation of the identity of
order $(\ez_1,\,\ez_2)$ (for short, continuous
$(\ez_1,\,\ez_2)$-$\ati$). For example, if $\{S_k\}_{k\in\zz}$ is an
$(\ez_1,\,\ez_2)$-$\ati$ and if we set $\wz S_t(x,\,y)=S_k(x,\,y)$
for $t\in(2^{-k-1},\,2^{-k}]$ with $k\in\zz$, then  $\{\wz
S_t\}_{t>0}$ is a continuous $(\ez_1,\,\ez_2)$-$\ati$.
\end{rem}

\begin{defn}\hspace{-0.2cm}{\bf.}\hspace{0.2cm}\label{d2.4}\rm
Let $\ez_1\in(0,\,1]$, $\ez_2\in(0, \fz)$
and $\{S_t\}_{t>0}$ be a continuous
$(\ez_1,\,\ez_2)-\ati$. Let $\rho$ be admissible.

(i) For any $f\in L^1_\loc(\cx)$
and $x\in\cx$, the radial maximal function $ S^+(f)$ is
defined by
\begin{equation*}
 S^+(f)(x)\equiv\sup_{t>0} |S_t(f)(x)|;
\end{equation*}

(ii) For any $f\in L^1_\loc(\cx)$ and $x\in\cx$, the radial maximal function $
S^+_{\rho}(f)$ associated to $\rho$ is defined by
\begin{equation*}
 S^+_{\rho}(f)(x)\equiv\sup_{0<t<\rho(x)} |S_t(f)(x)|.
\end{equation*}
\end{defn}

\vspace{-1cm}
\section*{}
\setcounter{section}{3}
\setcounter{thm}{0}
\setcounter{defn}{0}
\setcounter{lem}{0}
\setcounter{rem}{0}
\setcounter{prop}{0}
\setcounter{cor}{0}
\noindent 3. {\bf Localized BMO and BLO spaces.}
This section is divided into two subsections.
In Section 3.1, we introduce a localized BMO-type space $\bmo_\rho(\cx)$
and establish its several equivalent characterizations,
John-Nirenberg inequality and some other properties; while Section
3.2 is devoted to the study of a corresponding localized
BLO-type space $\blo_\rho(\cx)$.

\medskip
\noindent 3.1. {\bf A localized BMO space.}
\begin{defn}\hspace{-0.2cm}{\bf.}\hspace{0.2cm}\label{d3.1}\rm
 Let $\rho$ be an admissible function on $\cx$,
$\cd\equiv\{B(x,\,r)\subset\cx:\ x\in\cx,\
r\ge\rho(x)\}$ and $q\in[1,\,\fz)$.

(i)  A function $f\in
L^q_{\loc}(\cx)$ is said to be in the space $\bmo^q(\cx)$ if
$$\|f\|_{\bmo^q(\cx)}\equiv\sup_{B\subset\cx}
\lf\{\frac1{\mu(B)}\int_B|f(y)-f_B|^qd\mu(y)\r\}^{1/q}<\fz.$$

(ii) A function $f\in L^q_{\loc}(\cx)$ is said to be
in the space $\bmo^q_\rho(\cx)$ if
$$\begin{array}{cl}
\|f\|_{\bmo_\rho^q(\cx)}\equiv&\dsup_{B\notin\cd}
\lf\{\frac1{\mu(B)}\int_B|f(y)-f_B|^qd\mu(y)\r\}^{1/q}\\
&+\dsup_{B\in\cd}
\lf\{\frac1{\mu(B)}\int_B|f(y)|^qd\mu(y)\r\}^{1/q}<\fz.
\end{array}$$
\end{defn}

\begin{rem}\hspace{-0.2cm}{\bf.}\hspace{0.2cm}\label{r3.1}\rm
(i) The space $\bmo^q(\cx)$
with $q\in[1,\fz)$ coincides with $\bmo^1(\cx)$;
see \cite{cw77}. We denote $\bmo^1(\cx)$ simply
by $\bmo(\cx)$.

(ii) We also denote $\bmo^1_\rho(\cx)$ simply by $\bmoz$.
The localized space $\bmo_\rho(\rd)$
when $\rho\equiv1$ was first introduced by Goldberg \cite{g79}.
If $q>\frac d2$, $V\in \cb_q(\rd)$ and $\rho$ is as in \eqref{2.4},
then $\bmo_\rho(\rd)$ is just the space $\bmo_\cl(\rd)$ introduced by
 Dziub\'anski et al in \cite{dgmtz05}.
For all $q\in [1,\fz)$, $\bmo^q_\rho(\cx)\subsetneq \bmo(\cx)$.

(iii) Let $q\in[1, \fz)$,
$a\in(0, \fz)$ and $\cd_a\equiv\{B(x, r)\subset\cx: r\ge a\}$.
Define the space $\bmo^q_a(\cx)$ as in Definition \ref{d3.1} (ii) with
$\cd$ replaced by $\cd_a$.
Then, \eqref{2.1} implies that for all $q\in[1, \fz)$  and
fixed $a_1$, $a_2\in(0, \fz)$,
$\bmo^q_{a_1}(\cx)=\bmo^q_{a_2}(\cx)$ with equivalent norms.
From this, it further follows that
if $\mu(\cx)<\fz$, then for all $q\in[1, \fz)$ and any fixed
$a\in(0, \fz)$, $\bmo^q_\rho(\cx)=\bmo^q_a(\cx)$ with equivalent norms.
In fact, by \eqref{2.2},
there exists a positive constant $M$ such that for all
$x$, $y\in\cx$, $d(x, y)\le M$.
This together with Lemma 2.1 in \cite{yz08} implies that
there exist positive constants $C$ and $\wz C$ such that
for all $x\in \cx$, $C\le \rho(x)\le \wz C$.
Thus, for all $q\in[1, \fz)$,
$\bmo^q_{C}(\cx)\subset\bmo^q_\rho(\cx)\subset
\bmo^q_{\wz C}(\cx)$
which implies the desired conclusion.
\end{rem}

The following result follows from Definition \ref{d3.1}.

\begin{lem}\hspace{-0.2cm}{\bf.}\hspace{0.2cm}\label{l3.1}
Let $\rho$ be an admissible function on $\cx$
and $q\in[1,\fz)$. Then
$\bmo_\rho(\cx)=\bmo_\rho^q(\cx)$ with equivalent norms.
\end{lem}

\begin{pf}\rm
Assume that $f\in\bmo_\rho^q(\cx)$. Then by the H\"older inequality,
$f\in\bmo_\rho(\cx)$ and
$\|f\|_{\bmo_\rho(\cx)}\le\|f\|_{\bmo_\rho^q(\cx)}$. Conversely, if
$f\in\bmo_\rho(\cx)$, then by Remark
\ref{r3.1} (i) and Definition \ref{d3.1},
\begin{eqnarray*}
\|f\|_{\bmo_\rho^q(\cx)}\ls\|f\|_{\bmo^q(\cx)}+ \sup_{B\in\cd}|f_B|
&&\ls\|f\|_{\bmo_\rho(\cx)},
\end{eqnarray*}
which implies that $f\in\bmo_\rho^q(\cx)$ and
$\|f\|_{\bmo_\rho^q(\cx)}\ls\|f\|_{\bmo_\rho(\cx)}$. Thus
$\bmo_\rho(\cx)=\bmo_\rho^q(\cx)$ with equivalent
norms, which completes the proof of Lemma \ref{l3.1}.
\hfill$\Box$\end{pf}

Recall that the classical John-Nirenberg inequality (see
\cite{cw77}) says that there exist positive
constants $C_1$ and $C_2$ such that for all
$f\in\bmo(\cx)$, balls $B$ and $\lz>0$,
\begin{equation}\label{3.1}
\mu\lf(\{x\in B :\ |f(x)-f_B|>\lz\}\r)\le C_1\mu(B)
\exp\lf\{-\dfrac{C_2\lz}{\|f\|_{\bmo(\cx)}}\r\}.
\end{equation}
From this, we deduce a variant of the John-Nirenberg inequality
suitable for $\bmo_\rho(\cx)$ as follows.

\begin{thm}\hspace{-0.2cm}{\bf.}\hspace{0.2cm}\label{t3.1}
 Let $\rho$ be an admissible function on $\cx$
 and $\cd$ be as in Definition \ref{d3.1}.
 If $f\in\bmo_\rho(\cx)$, then there exist positive constants $C_3$
and $C_4$ such that for all balls $B$ and $\lz>0$,
\begin{equation}\label{3.2}
\mu\lf(\{x\in B:\ |f(x)-f_B|>\lz\}\r)\le C_3\mu(B)
\exp\lf\{-\dfrac{C_4\lz}{\|f\|_{\bmo_\rho(\cx)}}\r\},
\end{equation}
and, moreover, for all $B\in\cd$,
\begin{equation}\label{3.3}
\mu\lf(\{x\in B:\ |f(x)|>\lz\}\r)\le
C_3\mu(B)\exp\lf\{-\dfrac{C_4\lz}{\|f\|_{\bmo_\rho(\cx)}}\r\}.
\end{equation}
\end{thm}

\begin{pf}\rm
The inequality \eqref{3.2} follows from \eqref{3.1} and Definition \ref{d3.1}
directly. To show \eqref{3.3}, let $B\in\cd$. If
$\lz>2\|f\|_{\bmo_\rho(\cx)}$, by the definition,
we have $\lz>2|f|_B$. Thus
for all balls $B$ in $\cd$, we obtain
$$\mu\lf(\{x\in B :\ |f(x)|>\lz\}\r)
\le\mu\lf(\{x\in B : \ |f(x)-f_B|>\lz/2\}\r),$$
which together \eqref{3.2}
yields \eqref{3.3}; if $0<\lz\le2\|f\|_{\bmo_\rho(\cx)}$, we then
have
$$\mu\lf(\{x\in B :\ |f(x)|>\lz\}\r)\le\mu(B)\ls\mu(B)
\exp\lf\{-\dfrac{C_2\lz}{\|f\|_{\bmo_\rho(\cx)}}\r\}.$$ This
finishes the proof of Theorem \ref{t3.1}.
\hfill$\Box$\end{pf}

\begin{rem}\hspace{-0.2cm}{\bf.}\hspace{0.2cm}\label{r3.2}\rm
 Let $\rho$ be an admissible function on $\cx$ and $q\in[1,\fz)$.
Applying Theorem \ref{t3.1}, we can also obtain that $\bmo_\rho(\cx)=
\bmo_\rho^q(\cx)$ with equivalent norms.
\end{rem}

We now establish the relation between $\bmo(\cx)$ and $\bmoz$
in terms of certain approximation of the identity.
To begin with, let $\rho$ be an admissible function on $\cx$. In \cite{yz08},
it was proved that there exist a nonnegative function
$K_\rho$ on $\cx\times\cx$ and a positive constant $C_5$ such that

${\rm (K)_1}$ $K_\rho(x,\,y)=0$ if $x,\,y\in\cx$ satisfying
$d(x,\,y)>C_5\min\{\rho(x),\,\rho(y)\}$;

${\rm (K)_2}$ $K_\rho(x,\,y)\le C_5\frac1{V_{\rho(x)}(x)+V_{\rho(y)}(y)}$
for all $x,\,y\in\cx$;

${\rm (K)_3}$ $K_\rho(x,\,y)=K_\rho(y,\,x)$ for all $x,\,y\in\cx$;

${\rm (K)_4}$ $\int_\cx K_\rho(x,\,y)\,d\mu(x)=1$
for all $y\in\cx.$

For all $x\in\cx$, let
\begin{equation}\label{3.4}
K_\rho(f)(x)\equiv\int_\cx K_\rho(x,\,y)f(y)\,d\mu(y).
\end{equation}
It was proved in \cite{yz08} that if $f\in H^1_\rho (\cx)$,
the Hardy space associated to $\rho$, then $f-K_\rho (f)\in H^1(\cx)$,
where $H^1(\cx)$ is the Hardy space studied in \cite{hmy1,hmy2,gly1},
which coincides with the atomic Hardy space $H^1_{at}(\cx)$
of Coifman and Weiss in \cite{cw77}.
Moreover, there exists a positive constant $C$ such that
for all $f\in H^1_\rho (\cx)$,
\begin{equation}\label{3.5}
\|f-K_\rho (f)\|_{H^1(\cx)}\le C\|f\|_{H^1_\rho (\cx)}.
\end{equation}
On the other hand, it was showed in \cite{yyz08} that
the dual space of $H^1_\rho(\cx)$ is $\bmoz$. From these
facts, we deduce the following corollary.

\begin{cor}\hspace{-0.2cm}{\bf.}\hspace{0.2cm}\label{c3.1} Let $\rho$ be an
admissible function on $\cx$ and $K_\rho$ be as in \eqref{3.4}. Then

(i) $\bmoz=\{b\in\bmox:\ K_\rho(b)\in\lin\}$; moreover, for all
$b\in\bmoz$, $\|b\|_\bmoz\sim \|K_\rho(b)\|_\lin+\|b\|_\bmox$.

(ii) If $f\in\bmox$, then $f-K_\rho(f)\in\bmoz$; moreover, there exists
a positive constant $C$ such that for all $f\in\bmox$,
$$\|f-K_\rho(f)\|_\bmoz\le C\|f\|_\bmox.$$
\end{cor}

\begin{pf}\rm
We first prove (i). Assume that $b\in\bmox$ with $K_\rho(b)\in\lin$.
Recall that $H^1_\rho(\cx)\subset L^1(\cx)$ (see Lemma 3.1 in \cite{yz08}).
For any $f\in H^1_\rho(\cx)$, by ${\rm (K)_3}$, \eqref{3.5} and
$\lf(H^1(\cx)\r)^\ast=\bmox$ (see \cite{cw77}),
we have
\begin{eqnarray*}
\lf|\dint_\cx b(x)f(x)\,d\mu(x)\r|
&&\le\lf|\dint_\cx b(x)\lf[f(x)-K_\rho(f)(x)\r]\,d\mu(x)\r|+
\lf|\dint_\cx b(x)K_\rho(f)(x)\,d\mu(x)\r|\\
&&\le\|b\|_\bmox\|f-K_\rho(f)\|_{H^1(\cx)}+
\lf|\dint_\cx f(x)K_\rho(b)(x)\,d\mu(x)\r|\\
&&\ls\|f\|_{H^1_\rho(\cx)}\lf[\|b\|_\bmox+\|K_\rho(b)\|_\lin\r].
\end{eqnarray*}
Thus by $(H^1_\rho(\cx))^\ast=\bmoz$, we obtain $b\in\bmoz$ and
$$\|b\|_\bmoz\ls\|b\|_\bmox+\|K_\rho(b)\|_\lin.$$

Conversely, assume that $b\in\bmoz$. By ${\rm (K)_1}$
and ${\rm (K)_2}$, for all $x\in \cx$, we have
\begin{eqnarray*}
|K_\rho(b)(x)|&\le& \dint_{d(x, y)\le C_5\rho(x)}
|K_\rho(x, y)b(y)|\,d\mu(y)\\
&\ls&\frac1{V_{C_5\rho(x)}(x)}\dint_{B(x,\,C_5\rho(x))}|b(y)|\,d\mu(y)
\ls\|b\|_\bmoz.
\end{eqnarray*}
This shows (i).

To see (ii), by ${\rm (K)_3}$, \eqref{3.5} and $(H^1(\cx))^\ast=\bmox$,
we have that for all $f\in\bmox$ and $b\in H^1_\rho(\cx)$,
\begin{eqnarray*}
\lf|\dint_\cx[f(x)-K_\rho(f)(x)]b(x)\,d\mu(x)\r|
&&=\lf|\dint_\cx f(x)[b(x)-K_\rho(b)(x)]\,d\mu(x)\r|\\
&&\ls\|f\|_\bmox\|b\|_{H^1_\rho(\cx)},
\end{eqnarray*}
which together with $(H^1_\rho(\cx))^\ast=\bmoz$
 implies that $f-K_\rho(f)\in\bmoz$ and the desired estimate.
 This finishes the proof
 of Corollary \ref{c3.1}.
\hfill$\Box$\end{pf}

\medskip
\noindent 3.2. {\bf A localized BLO space.}

\begin{defn}\hspace{-0.2cm}{\bf.}\hspace{0.2cm}\label{d3.2}\rm
Let $\rho$ and $\cd$ be as in Definition \ref{d3.1} and $q\in[1,\fz)$.

(i) A function $f\in L^q_{\loc}(\cx)$ is said to be in the space
$\blo^q(\cx)$ if
$$\|f\|_{\blo^q(\cx)}\equiv
\sup_{B\subset \cx}\lf\{\frac1{\mu(B)}
\int_B\lf[f(y)-{\mathop\einf_B} f\r]^q\,d\mu(y)\r\}^{1/q}<\fz.$$

(ii) A function $f\in L^q_{\loc}(\cx)$ is said to be in the space
$\blo_\rho(\cx)$ if
 \begin{eqnarray*}
\|f\|_{\blo^q_\rho(\cx)}&&\equiv \sup_{B\notin
\cd} \lf\{\frac1{\mu(B)}
\int_B\lf[f(y)-{\mathop\einf_B} f\r]^q\,d\mu(y)\r\}^{1/q}\\
&&\hs+\sup_{B\in\cd}\lf\{\frac1{\mu(B)}
\int_B|f(y)|^q\,d\mu(y)\r\}^{1/q}<\fz.
\end{eqnarray*}
\end{defn}

\begin{rem}\hspace{-0.2cm}{\bf.}\hspace{0.2cm}\label{r3.3}\rm
(i) The space $\blo^1(\rd)$ with the Lebesgue measure
was introduced by Coifman and
Rochberg \cite{cr80}, and extended by Jiang \cite{j05}
to the setting of $\rd$ with a non-doubling measure.
Let $q\in[1, \fz)$. Then the facts that $\blo^1(\cx)\subset\bmo(\cx)=
\bmo^q(\cx)$ together with the H\"older inequality imply that
$\blo^q(\cx)=\blo^1(\cx)$ with equivalent norms.
We denote $\blo^1(\cx)$ simply by $\blo(\cx)$.
Notice that $\blo(\cx)$ is not a linear space.

(ii) We also denote $\blo^1_\rho(\cx)$ simply
by $\bloz$. The localized $\blo$ space was first introduced in
\cite{hyy} in the setting of $\rd$ with a non-doubling measure.
For all $q\in [1,\fz)$,
$\blo^q_\rho(\cx)\subset\bmo^q_\rho(\cx)$. Even when
$\rho\equiv1$, it is not so difficult to show that for
all $q\in [1, \fz)$, $\blo_\rho^q(\rd)$ is a proper subspace of $\bmo_\rho^q(\rd)$.

(iii) Let $q\in[1, \fz)$,
$a\in(0, \fz)$ and $\cd_a\equiv\{B(x, r)\subset\cx: r\ge a\}$.
Define the space $\blo^q_a(\cx)$ as in Definition \ref{d3.2} (ii) with
$\cd$ replaced by $\cd_a$.
If $\mu(\cx)<\fz$, then for all $q\in[1, \fz)$ and admissible functions $\rho$,
and any fixed
$a\in(0, \fz)$, $\blo^q_\rho(\cx)=\blo^q_a(\cx)$ with equivalent norms.
The proof is similar to that of Remark \ref{r3.1} (iii) and is omitted.
\end{rem}

The following result follows from Definitions \ref{d3.1} and
\ref{d3.2}, whose proof is similar to that of Lemma \ref{l3.1}
and is omitted.

\begin{lem}\hspace{-0.2cm}{\bf.}\hspace{0.2cm}\label{l3.2}
Let $\rho$ be an admissible function on $\cx$
and $q\in[1,\fz)$. Then
$\blo_\rho(\cx)=\blo_\rho^q(\cx)$ with equivalent norms.
\end{lem}

\begin{thm}\hspace{-0.2cm}{\bf.}\hspace{0.2cm}\label{t3.2}
 Let $\rho$ be an admissible function on $\cx$.
There exists a positive constant $C$ such that for all
$f\in \bloz$, $f(x)\ge -C\|f\|_\bloz$ for almost all $x\in\cx$.
Moreover, the following statements are
equivalent:

(i) $f\in \bloz$;

(ii) $f\in L^1_\loc(\cx)$ and
there exists a nonnegative constant $A$ such that
\begin{eqnarray*}
&&\sup_{B\notin\cd}\lf\{\frac1{\mu(B)}\int_B
\lf[f(y)-{\mathop\infi_{B}}f\r]d\mu(y)\r\}\\
&&\hs+\sup_{B\in\cd}\lf\{\frac1{\mu(B)}\int_B|f(y)|d\mu(y)+
\lf|{\mathop\einf_B}f\r|\r\}\le A;\nonumber
\end{eqnarray*}

(iii) $f\in L^1_\loc(\cx)$ and
there exists a nonnegative constant $\wz C$ such that
\begin{equation*}
\sup_{B\subset\cx}\lf\{\frac1{\mu(B)}\int_B
\lf[f(y)-{\mathop\infi_{B}}f\r]d\mu(y)\r\}
+\sup_{B\in\cd}\lf\{\frac1{\mu(B)}\int_B|f(y)|d\mu(y)\r\}\le \wz C.
\end{equation*}
Moreover, $\|f\|_\bloz$, $\inf\{A\}$ and $\inf\{\wz C\}$
are mutually equivalent.
\end{thm}

\begin{pf}\rm
Let $f\in \bloz$. For all balls $B\equiv B(x_0,\,\rho(x_0)/2)$, we have
that
$${\mathop\infi_{B} f} \ge -\frac1{\mu(B)}\int_{B}
\lf[f(x)-{\mathop\infi_{B} f}\r]\,d\mu(x)-
\frac1{\mu(B)}\int_{B}|f(x)|\,d\mu(x)\gs-\|f\|_{\bloz},$$ which
together with $\cx=\cup_xB(x,\,\rho(x)/2)$ and the Vitali-Wiener
type covering lemma (see \cite[p.\,623]{cw77})
implies that there exists certain positive constant $C$
such that for $\mu$-a.\,e. $x\in\cx$,  $f(x)\ge -C\|f\|_\bloz$.
From this, it is easy to see that (i)
implies (ii). Obviously, (ii) implies (iii) and (iii) implies (i).
Thus we complete the proof of Theorem \ref{t3.2}.
\hfill$\Box$\end{pf}

\begin{rem}\hspace{-0.2cm}{\bf.}\hspace{0.2cm}\rm\label{r3.4}
(i) From Theorem \ref{t3.2} (ii) and Definition \ref{d3.2} (i),
it follows that $\bloz\subset\blox$.

(ii) During this paper being written, we learnt that when
$V\in B_q(\rd)$ with $q>d/2$, and $\rho$ is as
in \eqref{2.4}, Theorem \ref{t3.2} (iii) was used, independently,
by Gao, Jiang and Tang \cite{gjt} to introduce the space
$\blo_\cl(\rd)$ corresponding to the Schr\"odinger operator
$\cl=-\Delta+V$.
\end{rem}

As a consequence of Theorem \ref{t3.2},
we have the following corollary.

\begin{cor}\hspace{-0.2cm}{\bf.}\hspace{0.2cm}\label{c3.2}
 Let $\rho$ be an admissible function on $\cx$. Then
$\bloz=\bmoz\cap\blox$. Moreover, for all $f\in \bloz$,
$$\|f\|_\bloz\sim \|f\|_\bmoz+\|f\|_\blox.$$
\end{cor}

\begin{pf}\rm
Let $f\in \bmoz\cap\blox$ first. By Definitions
\ref{d3.1} and \ref{d3.2}, $f\in \bloz$ and $\|f\|_\bloz\le
\|f\|_\bmoz+\|f\|_\blox$. Conversely, assume that $f\in \bloz$. It
follows from Definitions \ref{d3.1} and \ref{d3.2} that $f\in \bmoz$ and
$\|f\|_\bmoz\ls \|f\|_\bloz$, which together with Remark \ref{r3.4} (i)
completes the proof of Corollary \ref{3.2}.
\hfill$\Box$\end{pf}

As a consequence of Corollary \ref{c3.1} and Corollary \ref{c3.2}, we have the
following result.

\begin{cor}\hspace{-0.2cm}{\bf.}\hspace{0.2cm}\label{c3.3}
Let $\rho$ be an admissible function on $\cx$
and $K_\rho$ be as in \eqref{3.4}. Then
$$\bloz=\{f\in\blox:\ K_\rho(f)\in\lin\}.$$
\end{cor}

\begin{pf}\rm
Assume that $f\in \bloz$ first. Then by Corollary \ref{c3.2},
$f\in\blox\cap\bmoz$.
From this and Corollary \ref{c3.1} (i), it follows that $K_\rho(f)\in\lin$.
Conversely, if $f\in\blox$ and $K_\rho(f)\in\lin$,
then the obvious fact $\blox\subset\bmox$ together with another application of Corollary \ref{c3.1} (i)
implies that $f\in\bmoz$, which together with Corollary
\ref{c3.2} yields that $f\in\bloz$.
This finishes the proof of Corollary \ref{c3.3}.
\hfill$\Box$\end{pf}

\begin{thm}\hspace{-0.2cm}{\bf.}\hspace{0.2cm}\label{t3.3}
Let $\rho$ be an admissible function on $\cx$
and $K_\rho$ be as in \eqref{3.4}. Then there exists a
positive constant $C$ such that for all
$f\in\blo(\cx)$, $f-K_\rho f\in\bloz$ and
$$\|f-K_\rho(f)\|_\bloz\le C\|f\|_{\blo(\cx)}.$$
\end{thm}

\begin{pf}\rm
Let $f\in\blo(\cx)$.
By the homogeneity of $\|\cdot\|_\blox$ and $\|\cdot\|_\bloz$, we may assume that
$\|f\|_{\blo(\cx)}=1$. Let $B\equiv B(x_0, r)\in\cd$.
Observe that by \eqref{2.3}, for any $a\in(0, \fz)$,
there exists a constant $\wz C_a\in[1, \fz)$ such that
for all $x$, $y\in\cx$ with $d(x, y)\le a\rho(x)$,
\begin{equation}\label{3.7}
\rho(y)/\wz {C}_a\le\rho(x)\le {\wz C}_a\rho(y).
\end{equation}
By this and $r\ge \rho(x_0)$, we obtain that
for all $x\in B$, $\rho(x)\ls r$.
Then there exists a positive constant $C$ such that
for all $x\in B$,
$B(x, C_5\rho(x))\subset CB$. By ${\rm (K)_1}$
through ${\rm (K)_4}$, \eqref{2.1} and the Tonelli theorem,
we obtain
\begin{eqnarray}\label{3.8}
&&\frac1{\mu(B)}\dint_B\lf|f(x)-K_\rho(f)(x)\r|\,d\mu(x)\nonumber\\
&&\hs\le\frac1{\mu(B)}\dint_B\lf\{\lf[f(x)-{\mathop\einf_{CB}}f\r]
+\lf|K_\rho\lf(f-{\mathop\einf_{CB}}f\r)(x)\r|\r\}\,d\mu(x)\nonumber\\
&&\hs\ls1+\frac1{\mu(B)}\dint_{CB}\dint_{B(y,\,C_5\rho(y))}
\frac{|f(y)-{\mathop\einf_{CB}}f|}{V_{\rho(x)}(x)+V_{\rho(y)}(y)}
\,d\mu(x)\,d\mu(y)\ls1.
\end{eqnarray}

On the other hand, let $B\equiv B(x_0, r)\notin\cd$.
Using $r<\rho(x_0)$ and \eqref{3.7} with $a=1$, we obtain that there
exists a constant $\wz A_1\in [1,\fz)$
such that for all $x$, $y\in B$, $B(x, \rho(x))\subset B(y, \wz A_1\rho(y))$.
From this together with \eqref{2.1}, it follows that for all $x,$ $y\in B$,
$$\lf|{\mathop\einf_{B(y,\,\rho(y))}}f
-{\mathop\einf_{B(x,\,\rho(x))}}f\r|\ls1.$$
By this together with ${\rm (K)}_1$ through
$\rm (K)_4$ and \eqref{2.1}, we have that for all $x,$ $y\in B$,
\begin{eqnarray*}
\lf|K_\rho(f)(y)-{\mathop\einf_{B(x,\,\rho(x))}}f\r|
\le\lf|K_\rho(f)(y)-{\mathop\einf_{B(y,\,\rho(y))}}f\r|
+\lf|{\mathop\einf_{B(y,\,\rho(y))}}f
-{\mathop\einf_{B(x,\,\rho(x))}}f\r|\ls 1.
\end{eqnarray*}
From this fact, we deduce that
\begin{eqnarray*}
&&\lf[f-K_\rho(f)\r]_B-{\mathop\einf_{B}}\lf[f-K_\rho(f)\r]\\
&&\hs\le \frac1{\mu(B)}\dint_B\lf\{[f(x)-K_\rho(f)(x)]
-{\mathop\einf_{B}}f-{\mathop\einf_{B}}[-K_\rho(f)]\r\}\,d\mu(x)\\
&&\hs\ls 1+\frac1{\mu(B)}\dint_B\lf\{\lf[-K_\rho(f)(x)
+{\mathop\einf_{B(x,\,\rho(x))}}f\r]+\lf[-{\mathop\einf_{B(x,\,\rho(x))}}f
-{\mathop\einf_{B}}[-K_\rho(f)]\r]\r\}\,d\mu(x)\\
&&\hs\ls 1.
\end{eqnarray*}
This together with \eqref{3.8} gives the desired estimate
and hence, finishes the proof of Theorem \ref{t3.3}.
\hfill$\Box$\end{pf}

\vspace{-1.4cm}
\section*{}
\setcounter{section}{4}
\setcounter{thm}{0}
\setcounter{defn}{0}
\setcounter{lem}{0}
\setcounter{rem}{0}
\setcounter{prop}{0}
\setcounter{cor}{0}
\noindent 4. {\bf Boundedness of the natural and the Hardy-Littlewood
maximal functions.}
In this section, we first obtain the boundedness of
the natural maximal function, the Hardy-Littlewood
maximal function and their localized versions from
$\bmoz$ to $\bloz$; as an application, we then establish
several equivalent characterizations for $\blo_\rho(\cx)$
via the localized natural maximal function.

\begin{defn}\hspace{-0.2cm}{\bf.}\hspace{0.2cm}\label{d4.1}\rm
Let $\rho$ be an admissible function on $\cx$. For all $f\in L^1_\loc(\cx)$
and $x\in\cx$, define
$\hl(f)(x)\equiv\sup_{x\in B}\{|f|_B\}$, $M(f)(x)\equiv\sup_{x\in B}\{f_{B}\}$,
$\hl_{\rho}(x)\equiv\sup_{x\in B,\,B\notin\cd}\{|f|_{B}\}$, and
$M_{\rho}(f)(x)\equiv\sup_{x\in B,\,B\notin\cd}\{f_{B}\}.$
\end{defn}

\begin{thm}\hspace{-0.2cm}{\bf.}\hspace{0.15cm}\label{t4.1}
Let $\rho$ be an admissible function on $\cx$. Then
$M_{\rho}$ is bounded from $\bmo_{\rho}(\cx)$ to $\blo_{\rho}(\cx)$,
namely, there exists a positive
constant $C$ such that
for all $f\in \bmo_{\rho}(\cx)$,
$M_{\rho}(f)\in\blo_{\rho}(\cx)$ and
\begin{eqnarray*}
\|M_{\rho}(f)\|_{\blo_{\rho}(\cx)} \le C\|f\|_{\bmo_{\rho}(\cx)}.
\end{eqnarray*}
\end{thm}

\begin{pf}\rm
Let $f\in\bmoz$. By the homogeneity of $\|\cdot\|_\bmoz$ and
$\|\cdot\|_\bloz$, we may assume that $\|f\|_\bmoz=1$. We first
prove that for all balls $B\equiv B(x_0,\,r)\in\cd$,
\begin{equation}\label{4.1}
[\hl_{\rho}(f)]_B\ls1.
\end{equation}
From this, it follows that for all balls $B\in\cd$,
$(|M_{\rho}(f)|)_B\ls1$, which implies that
$|M_{\rho}(f)(x)|<\fz$
for $\mu$-a.\,e. $x\in\cx$.

To prove \eqref{4.1}, for all balls $B\in\cd$, write
$$[\hl_{\rho}(f)]_B\le
[\hl_\rho(f\chi_{3B})]_B +[\hl_{\rho}(f\chi_{(3B)^\complement})]_B.$$ The
H\"older inequality together with the $L^2(\cx)$-boundedness of
$\hl$ (see \cite{cw71}) and Lemma \ref{l3.1} gives us that
\begin{equation}\label{4.2}
[\hl_\rho(f\chi_{3B})]_B\le \{([\hl(f\chi_{3B})]^2)_B\}^{1/2} \ls
\{(|f\chi_{3B}|^2)_B\}^{1/2}\ls1.
\end{equation}

We now claim that for all $y\in B$,
$\hl_{\rho}(f\chi_{(3B)^\complement})(y)\ls 1$. In fact, for all
balls $\wz B\equiv B(\wz x,\,\wz r)\notin\cd$ such that $\wz B\cap
B\ne\emptyset$, we have that either $\wz B\subset 3B$ or $B\subset
3\wz B$. If $\wz B\subset 3B$, then
$(|f\chi_{(3B)^\complement}|)_{\wz B}=0$. If $B\subset 3\wz B$, then
by \eqref{3.7}, there exists a constant $\wz C\in[1,\,\fz)$
 such that $\wz C\wz r\ge\rho(\wz x)$.
Thus, using \eqref{2.1}, we have that $(|f\chi_{(3B)^\complement}|)_{\wz
B}\ls (|f|)_{\wz C\wz B}\ls1.$ The claim then follows from the two
estimates above, which together with \eqref{4.2} leads to
\eqref{4.1}.

We now prove that there exists a positive constant $C$ such that for
all balls $B\equiv B(x_0,r)\notin\cd$,
\begin{equation}\label{4.3}
\lf[M_{\rho}(f)\r]_B\le C+ {\mathop\einf_B} M_{\rho}(f).
\end{equation}
Write
$$[M_{\rho}(f)]_B\le \lf\{M_{\rho}\lf[(f-f_{B})\chi_{3B}\r]\r\}_B+
\lf\{M_\rho\lf[f_{B}\chi_{3B}+f\chi_{(3B)^\complement}\r]\r\}_B.$$
Using the H\"older inequality, the $L^2(\cx)$-boundedness of $\hl$,
\eqref{2.1} and Lemma \ref{l3.1}, we obtain that
\begin{equation}\label{4.4}
\{M_{\rho}[(f-f_{B})\chi_{3B}]\}_B \le \{\hl[(f-f_{B})\chi_{3B}]\}_B
\ls[(|f-f_{3B}|^2)_{3B}]^{1/2}\ls1.
\end{equation}

Now we show that for all $x$, $y\in B$,
$M_\rho(f_{B}\chi_{3B}+f\chi_{(3B)^\complement})(x) \le \wz
C+M_{\rho}(f)(y).$ For all balls $\wz B\equiv B(\wz x,\,\wz
r)\notin\cd$ containing $x$, we have either $\wz B\subset 3B$ or
$B\subset 3\wz B$. Assume that $\wz B\subset 3B$ first. Then
\begin{equation}\label{4.5}
\lf[f_{B}\chi_{3B}+f\chi_{(3B)^\complement}\r]_{\wz B} =f_{B}\le
M_{\rho}(f)(y).
\end{equation}
Now we assume that $B\subset 3\wz B$. If $3\wz B\notin\cd$, then the
fact $y\in B\subset 3\wz B$ gives us that $f_{3\wz B}\le M_{\rho}(f)(y)$,
which together with \eqref{2.1} implies that
\begin{eqnarray*}
&&\lf[f_{B}\chi_{3B}+f\chi_{(3B)^\complement}\r]_{\wz B}-M_{\rho}(f)(y)\\
&&\hs=\lf[\lf(f_B-f_{3\wz B}\r)\chi_{3B}+
\lf(f-f_{3\wz B}\r)\chi_{(3B)^\complement}\r]_{\wz B}+f_{3\wz B}- M_{\rho}(f)(y)\\
&&\hs\le \lf[|f_B-f_{3\wz B}|\chi_{3B}\r]_{\wz B} +\lf(|f-f_{3\wz
B}|\r)_{\wz B}\ls1.
\end{eqnarray*}
If $3\wz B\in\cd$, then $\rho(\wz x)\le 3\wz r$. Since $B$, $\wz B\notin\cd$,
$x,\,y\in B$ and $\wz x$, $x\in \wz B$, by \eqref{3.7},
$\rho(y)\sim\rho(x_0)\sim \rho(x)\sim\rho(\wz x)>\wz r$, which
implies that $\rho(y)\sim\rho(\wz x)\sim \wz r$. Let $\wz A\in [1, \fz)$ satisfying
$\rho(y)>\wz r/\wz A$. Then $B(y,\,\wz r/\wz A)\notin \cd$. By
the fact that $f_{B(y,\,\wz r/\wz A)}\le M_\rho(f)(y)$ together with
\eqref{2.1}, we have that
\begin{eqnarray*}
&&\lf[f_{B}\chi_{3B}+f\chi_{(3B)^\complement}\r]_{\wz B}-M_{\rho}(f)(y)\\
&&\hs\le\lf[\lf(f_B-f_{3\wz B}\r)\chi_{3B}+ \lf(f-f_{3\wz
B}\r)\chi_{(3B)^\complement}\r]_{\wz B} +\lf|f_{3\wz B}-f_{B(y,\,\wz
r/\wz A)}\r| \ls 1.
\end{eqnarray*}
Combining the two inequalities above and \eqref{4.5} leads to that
$$[M_\rho(f_{B}\chi_{3B}+f\chi_{(3B)^\complement})]_B
-\einf_BM_{\rho}(f)\ls 1,$$
which together with \eqref{4.4} further
implies  \eqref{4.3}. This finishes the proof of Theorem
\ref{t4.1}.
\hfill$\Box$\end{pf}

\begin{lem}\hspace{-0.2cm}{\bf.}\hspace{0.2cm}\label{l4.1}
Let $\rho$ be an admissible function on $\cx$.
Then $f\in\bloz$ if and only if $f\in L^1_\loc(\cx)$,
$M_\rho(f)-f\in\lin$ and
\begin{equation}\label{4.6}
\sup_{B\in\cd}\frac1{\mu(B)}\dint_B|f(y)|\,d\mu(y)\ls 1.
\end{equation}
Furthermore, $\|M_\rho(f)-f\|_\lin\sim\|f\|_\bloz.$
\end{lem}

\begin{pf}\rm
Assuming that $f\in\bloz$, we then see that \eqref{4.6} holds.
Since $\mu$ is regular, for $\mu$-a.\,e.\, $x\in\cx$, there exists a
sequence of balls $\{B_k\}_k$ centered at $x$ with $r_{B_k}\to 0$ as
$k\to \fz$ such that
\begin{equation}\label{4.7}
\lim_{k\to \fz}\frac1{\mu(B_k)}\dint_{B_k}f(y)\,d\mu(y)=f(x).
\end{equation}
Let $x$ be any point satisfying \eqref{4.7} and $B$ be a ball
containing $x$ with $B\notin \cd$. Then we obtain that $f(x)\ge
\mathop\einf_Bf$ and $f_B-f(x)\le\|f\|_\bloz.$ Taking the
supremum over all balls $B$ containing $x$ and $B\notin\cd$, we have
$M_\rho(f)(x)-f(x)\le \|f\|_\bloz.$

Conversely, assume that $f$ satisfies \eqref{4.6} and
$M_\rho(f)-f\in\lin$. Then for all balls $B\not\in\cd$ and
$\mu$-a.\,e. $x\in B$, $f(x)\ge f_B-\|M_\rho(f)-f\|_\lin.$ This
yields that
$$\mathop\einf_B f\ge f_B-\|M_\rho(f)-f\|_\lin,$$
which together with \eqref{4.6} implies that $f\in\bloz$ and
$\|f\|_\bloz\ls \|M_\rho(f)-f\|_\lin.$ This finishes the proof of
Lemma \ref{l4.1}.
\hfill$\Box$\end{pf}

\begin{thm}\hspace{-0.2cm}{\bf.}\hspace{0.2cm}\label{t4.2}
Let $\rho$ be an admissible function on $\cx$. Then
$f\in\bloz$ if and only if there exist $h\in\lin$ and
$g\in \bmoz$ such that
\begin{equation}\label{4.8}
f=M_\rho(g)+h.
\end{equation}
 Furthermore,
$\|f\|_\bloz\sim \inf(\|g\|_\bmoz+\|h\|_\lin),$ where the infimum is
taken over all representations of $f$ as in \eqref{4.8}.
\end{thm}

\begin{pf}\rm
If there exist $g$ and $h$ satisfying \eqref{4.8}, then by Theorem
\ref{t4.1}, $M_\rho(g)\in\bloz$, which together with
$L^\fz(\cx)\subset\bloz$ implies that $f\in\bloz$ and
$$\|f\|_\bloz\ls \|M_\rho(g)\|_\bloz+\|h\|_\lin\ls\|g\|_\bmoz+\|h\|_\lin.$$

To see the converse, assume that $f\in\bloz$. By $\bloz\subset\bmoz$ and
Theorem \ref{t4.1}, we see $M_\rho(f)\in\bloz$. Let $h\equiv f-M_\rho(f)$ and
$g\equiv f$. Then Theorem \ref{t4.2} follows from Lemma \ref{l4.1},
which completes the proof of Theorem \ref{t4.2}.
\hfill$\Box$\end{pf}

As another corollary of Theorem \ref{t4.1}, we obtain the
boundedness of $\hl,\,\hl_{\rho}$ and $M$ from $\bmo_{\rho}(\cx)$ to
$\blo_{\rho}(\cx)$. To this end, we first establish the following
useful lemma.

\begin{lem}\hspace{-0.2cm}{\bf.}\hspace{0.2cm}\label{l4.2}
Let $\rho$ be an admissible function on $\cx$ and
 $Y$ be one of the spaces $\bmo(\cx)$, $\bmo_{\rho}(\cx)$,
$\blo(\cx)$ and $\blo_{\rho}(\cx)$. If $f\in Y$,  $h_1,\, h_2\in
L^\fz(\cx)$, and $f-h_2\le g\le f+h_1$, then $g\in Y$ and
$$\|g\|_Y\le \|f\|_Y+\|h_1\|_{L^\fz(\cx)}+\|h_2\|_{L^\fz(\cx)}.$$
\end{lem}

\begin{pf}\rm
We only consider the case that $Y=\blo_{\rho}(\cx)$ by similarity.
For all balls $B\in\cd$, we have that
\begin{eqnarray*}
|g|_B\le |f|_B+\|h_1\|_{L^\fz(\cx)}+\|h_2\|_{L^\fz(\cx)} \le
\|f\|_{\blo_{\rho}} +\|h_1\|_{L^\fz(\cx)}+\|h_2\|_{L^\fz(\cx)}.
\end{eqnarray*}
On the other hand, for all balls $B\notin \cd$,
$$g_B-{\mathop\einf_B}g\le f_B+\|h_1\|_{L^\fz(\cx)}
-{\mathop\einf_B}f+\|h_2\|_{L^\fz(\cx)} \le \|f\|_{\blo_{\rho}}+
\|h_1\|_{L^\fz(\cx)}+\|h_2\|_{L^\fz(\cx)}.$$
 Combining the two
inequalities above finishes the proof of Lemma \ref{l4.2}.
\hfill$\Box$\end{pf}

\begin{cor}\hspace{-0.2cm}{\bf.}\hspace{0.2cm}\label{c4.1}
Let $\rho$ be an admissible function on $\cx$. Then $\hl,\,\hl_{\rho}$ and
$M$ are bounded from $\bmo_{\rho}(\cx)$ to $\blo_{\rho}(\cx)$,
namely, there exists a positive constant $C$ such that for all $f\in
\bmo_{\rho}(\cx)$,
$\hl(f),\,\hl_{\rho}(f),\,M(f)\in\blo_{\rho}(\cx)$ and
\begin{eqnarray*}
\|\hl(f)\|_{\blo_{\rho}(\cx)}+\|\hl_{\rho}(f)\|_{\blo_{\rho}(\cx)}+
\|M(f)\|_{\blo_{\rho}(\cx)}\le C\|f\|_{\bmo_{\rho}(\cx)}.
\end{eqnarray*}
\end{cor}

\begin{pf}\rm
Since for all locally integrable functions $f$,
$$M_{\rho}(f)\le M(f)\le M_{\rho}(f)
+\|f\|_{\bmo_{\rho}(\cx)},$$ by Theorem \ref{t4.1} and Lemma
\ref{l4.2}, we have that if $f\in\bmoz$, then $M(f)\in\bloz$ and
$\|M(f)\|_{\blo_{\rho}(\cx)}\ls \|f\|_{\bmo_{\rho}(\cx)}.$ Using
this together with Theorem \ref{t4.1}, the facts that
$\hl_{\rho}(f)=M_{\rho}(|f|)$, $\hl(f)=M(|f|)$ and that
$\||f|\|_{\bmo_{\rho}(\cx)} \le 2\|f\|_{\bmo_{\rho}(\cx)}$ for all
$f\in\bmoz$, we have that
$\hl(f),\,\hl_{\rho}(f)\in\blo_{\rho}(\cx)$ and
$$\|\hl_{\rho}(f)\|_{\blo_{\rho}(\cx)}
+\|\hl(f)\|_{\blo_{\rho}(\cx)} \ls \|f\|_{\bmo_{\rho}(\cx)}.$$
This finishes the proof of Corollary \ref{c4.1}.
\hfill$\Box$\end{pf}

\begin{rem}\hspace{-0.2cm}{\bf.}\hspace{0.2cm}\label{r4.1} \rm
Let $d\ge 3$, $V$ be a nonnegative integrable function on $\rd$
and $\cl=-\Delta+V$.
If $q>d/2$, $V\in\cb_q(\rd)$ and $\rho$ is
as in \eqref{2.4}, then $\bmo_\rho(\rd)$ is just the
space $\bmo_\cl(\rd)$  introduced in
 \cite{dgmtz05}. It was proved in \cite{dgmtz05} that
$\hl$ is bounded on $\bmo_\rho(\rd)$.
Recall that $\blo_\rho(\rd)\subsetneq\bmo_\rho(\rd)$.
Thus, Corollary \ref{c4.1} improves the
result of \cite{dgmtz05}. Similar claim is also true
for $\hl$ on Heisenberg groups; see \cite{ll08}.
\end{rem}

\vspace{-1.4cm}
\section*{}
\setcounter{section}{5}
\setcounter{thm}{0}
\setcounter{defn}{0}
\setcounter{lem}{0}
\setcounter{rem}{0}
\setcounter{prop}{0}
\setcounter{cor}{0}
\noindent 5. {\bf Boundedness of several maximal operators.}
This section consists of two subsections.
Subsection 5.1 is devoted to the boundedness
of several radial maximal operators from $\bmo_\rho(\cx)$
to $\bloz$; while in Subsection 5.2, we
obtain the boundedness of the Poisson semigroup maximal operator
from $\bmoz$ to $\bloz$.

\medskip
\noindent 5.1. {\bf Boundedness of radial maximal operators.}
\begin{thm}\hspace{-0.2cm}{\bf.}\hspace{0.2cm}\label{t5.1}
Let $\rho$ be an admissible function on $\cx$ and $S^+(f)$ be as in
Definition \ref{d2.4}. Then there exists a positive constant $C$ such
that for all $f\in\bmo_{\rho}(\cx)$,
$$\|S^+(f)\|_{\blo_{\rho}(\cx)}\le C\|f\|_{\bmo_{\rho}(\cx)}.$$
\end{thm}

\begin{pf}\rm
By the homogeneity of $\|\cdot\|_\bmoz$ and $\|\cdot\|_\bloz$, we
may assume that
$$\|f\|_{\bmo_{\rho}(\cx)}=1.$$
Observe that by
Definition \ref{d2.3} (i), $S^+(f)\ls \hl(f)$. From this and
Corollary \ref{c4.1}, it follows that for all balls $B\equiv B(x_0,\,r)\in\cd$,
$$\frac1{\mu(B)}\int_B S^+(f)(x)\,d\mu(x)\ls1.$$
This also implies that $S^+(f)(x)<\fz$ for $\mu$-a.\,e. $x\in\cx$.
Moreover, by the inequality above, to finish the proof Theorem
\ref{t5.1}, it suffices to show that for all balls $B\equiv
B(x_0,\,r)\notin\cd$ and $y\in B$,
\begin{equation}\label{5.1}
\frac1{\mu(B)}\int_B[S^+(f)(x)-S^+(f)(y)]\,d\mu(x)\ls1.
\end{equation}

Let $f_1\equiv(f-f_B)\chi_{2B}$,
$f_2\equiv(f-f_B)\chi_{(2B)^\complement},$ $B_1\equiv\lf\{x\in B:\
S_{r}^+(f)(x)\ge S_{\fz}^+(f)(x)\r\}$ and $B_2\equiv B\setminus
B_1$, where for all $x\in\cx$,
$S_{r}^+(f)(x)\equiv\sup_{0<t<r}|S_t(f)(x)|$ and
$S_{\fz}^+(f)(x)\equiv\sup_{r\le t<\fz}|S_t(f)(x)|.$
By using Definition \ref{d2.3} (iv), we have that for
all $y\in B$,
\begin{eqnarray*}
&&\frac1{\mu(B)}\int_B[S^+(f)(x)-S^+(f)(y)]\,d\mu(x)\\
&&\hs=\frac1{\mu(B)}\int_{B_1}\lf[S_{r}^+(f)(x)-S^+(f)(y)\r]\,d\mu(x)\\
&&\hs\hs
+\frac1{\mu(B)}\int_{B_2}\lf[S_{\fz}^+(f)(x)-S^+(f)(y)\r]\,d\mu(x)\\
&&\hs\le \frac1{\mu(B)}\int_BS_{r}^+(f_1)(x)\,d\mu(x)+
\frac1{\mu(B)}\int_BS_{r}^+(f_2)(x)\,d\mu(x)\\
&&\hs\hs+|f_B-S_r(f)(y)|
+\frac1{\mu(B)}\int_{B_2}\lf[S_{\fz}^+(f)(x)
-S_{\fz}^+(f)(y)\r]\,d\mu(x)\\
&&\hs\equiv{\rm L_1}+{\rm L_2}+{\rm L_3}+{\rm L_4}.
\end{eqnarray*}

By the H\"older inequality, $S^+(f)\ls \hl(f)$, the
$L^2(\cx)$-boundedness of $\hl$, \eqref{2.1} and Lemma \ref{l3.1},
we obtain
$${\rm L_1}\ls\lf( \frac1{\mu(B)}\int_B[\hl(f_1)(x)]^2\,d\mu(x)\r)^{1/2}\ls
\lf( \frac1{\mu(B)}\int_{2B}|f(x)-f_{2B}|^2\,d\mu(x)\r)^{1/2}\ls1.$$

Recall that $\{S_t\}_{t>0}$ is a continuous $(\ez_1,\,\ez_2)$-$\ati$.
By Remark \ref{r2.1}, Definition \ref{d2.3} (i), \eqref{2.1} and the
fact that for all $x\in B$ and $j\in\nn$, $2^{j+1}B\subset B(x, 2^{j+2}r)$, we have that
for all $t\in(0, r)$,
\begin{eqnarray*}
|S_t(f_2)(x)|&&\ls\dint_{(2B)^\complement}\frac1{V(x, z)}\lf(\frac{t}{t+d(x, z)}\r)^{\ez_2}
|f(z)-f_B|\,d\mu(z)\\
&&\ls\sum_{j=1}^\fz \lf(\frac{t}{2^{j-1}r}\r)^{\ez_2}
\frac1{V_{2^{j-1}r}(x)}\int_{2^{j+1}B}|f(z)-f_B|\,d\mu(z)\\
&&\ls\lf(\frac{t}{r}\r)^{\ez_2}\sum_{j=1}^\fz 2^{-j\ez_2}
\lf[\frac1{\mu(2^{j+1}B)}\int_{2^{j+1}B}
|f(z)-f_{2^{j+1}B}|\,d\mu(z)
+|f_B-f_{2^{j+1}B}|\r]\\
&&\ls\lf(\frac{t}{r}\r)^{\ez_2}\sum_{j=1}^\fz 2^{-j\ez_2}(j+1)\ls1,
\end{eqnarray*}
which implies that ${\rm L_2}\ls1$.

By Definition \ref{d2.3} (iv) and (i) together with \eqref{2.1}
and the fact that for all $y\in B$ and $j\in\nn\cup\{0\}$,
$2^{j+1}B\subset B(y,\,2^{j+2}r)$, we have
\begin{eqnarray*}
{\rm L_3}&&\le\int_\cx \lf|S_r(y,\,z)[f(z)-f_B]\r|\,d\mu(z)\\
&&\ls\int_\cx\frac1{V_r(y)+V(y,\,z)}
\lf(\frac r{r+d(y,\,z)}\r)^{\ez_2}|f(z)-f_B|\,d\mu(z)\\
&&\ls\sum_{j=0}^\fz 2^{-j\ez_2}
\frac1{V_{2^{j-1}r}(y)}\int_{2^{j+1}B}|f(z)-f_B|\,d\mu(z)
\ls\sum_{j=0}^\fz 2^{-j\ez_2}(1+j)\ls1.
\end{eqnarray*}

On the other hand, for all $x,\,y\in B$ and $t\in[r, \fz)$,
$B(x,\,t)\subset B(y,\,2t)\subset B(x,\,3t)$. It follows from this
fact and \eqref{2.1} that
\begin{equation*}
\lf|f_{B(x,\,t)}-f_{B(y,\,t)}\r|\ls \frac1{\mu(B(x,\,3t))}
\dint_{B(y,\,2t)}\lf|f(z)-f_{B(y,\,t)}\r|\,d\mu(z)\ls1.
\end{equation*}
By this and an argument similar to the estimate for ${\rm L_3}$, we have
that for all $x$, $y\in B$ and $t\in[r, \fz)$,
\begin{eqnarray*}
&&|S_t(f)(x)-S_t(f)(y)|\\
&&\hs\ls \lf|S_t(f)(x)-f_{B(x,\,t)}\r|+
\lf|f_{B(x,\,t)}-f_{B(y,\,t)}\r|+\lf|f_{B(y,\,t)}-S_t(f)(y)\r|\ls1,
\end{eqnarray*}
which implies that
$${\rm L_4}\ls \frac1{\mu(B)}\int_B\sup_{r\le t<\fz}|S_t(f)(x)-S_t(f)(y)|\,d\mu(x)
\ls1.$$ Combining the estimates for ${\rm L_1}$ through ${\rm L_4}$
yields \eqref{5.1},
which completes the proof of Theorem \ref{t5.1}.
\hfill$\Box$\end{pf}

By Definition \ref{d2.3} (i), we have that for all $x\in \cx$ and
$t\in[\rho(x), \fz)$,
\begin{equation}\label{5.2}
\begin{array}[b]{ccl}
\lf|S_t(f)(x)\r|&&\ls\dint_\cx\frac1{V_t(x)+V(x,\,y)}
\lf(\dfrac{t}{t+d(x,\,y)}\r)^{\ez_2}|f(y)|\,d\mu(y)\\
&&\ls\dsum_{j=0}^\fz2^{-j\ez_2}\dfrac1{V_{2^{j-1}t}(x)}
\int_{d(x,\,y)<2^jt} |f(y)|\,d\mu(y)\ls\|f\|_\bmoz.
\end{array}
\end{equation}
This implies that there exists a positive constant $\wz C$ such that
for all $f\in\bmoz$,
$$S^+_{\rho}(f)\le S^+(f)\le S^+_{\rho}(f)+\wz C \|f\|_\bmoz.$$
From this, Lemma \ref{l4.2} and Theorem \ref{t5.1},
we deduce the following corollary.

\begin{cor}\hspace{-0.2cm}{\bf.}\hspace{0.2cm}\label{c5.1}
Let $\rho$ be an admissible function on $\cx$ and $S^+_\rho(f)$ be
as in Definition \ref{d2.4}. Then there exists a positive constant $C$
such that for all $f\in\bmo_{\rho}(\cx)$, $S^+_\rho(f)\in\bloz$ and
$\|S^+_\rho(f)\|_{\blo_{\rho}(\cx)}\le C\|f\|_{\bmo_{\rho}(\cx)}.$
\end{cor}

Let $\{T_t\}_{t>0}$ be a family of bounded linear
  operators with integral kernels $\{T_t(x,\,y)\}_{t>0}$.
Assume that there exist constants $C\in(0,\,\fz)$,  $\ez_1\in(0,1]$,
$\ez_2\in(0,\,\fz)$, $\dz\in(0,\,1]$ and $\gz\in(0,\,\fz)$, and an
$(\ez_1, \ez_2)$-$\ati$ $\{\wz T_t\}_{t>0}$ with kernels $\{\wz T_t(x,
y)\}_{t>0}$ such that for all $t\in(0, \fz)$ and $x,\,
y\in\cx$,
\begin{equation}\label{5.3}
\lf|T_t(x,\,y)-\wz T_t(x,\,y)\r|\le
C\lf(\frac{t}{t+\rho(x)}\r)^{\dz} \frac1{V_t(x)+V(x,\,y)}\lf(\frac
t{t+d(x,\,y)}\r)^{\gz}.
\end{equation}

Notice that by \eqref{5.3}, for all $f\in\bmoz$ and
$t\in(0,\,\fz)$, we have that for all $x\in\cx$,
\begin{eqnarray}\label{5.4}
\lf|T_t(f)(x)-\wz T_t(f)(x)\r|&&\ls \lf(\frac{t}{t+\rho(x)}\r)^{\dz}
\sum_{j=1}^\fz2^{-j\gz}\frac1{V_{2^{j-1}t}(x)}\int_{d(x,\,y)<2^jt}
 f(y)\,d\mu(y)\nonumber\\
&&\ls \|f\|_\bmoz\lf(\frac{t}{t+\rho(x)}\r)^{\dz}
\sum_{j=1}^\fz2^{-j\gz}\log\lf(1+\frac {\rho(x)}{2^jt}\r)\nonumber\\
&&\ls \|f\|_\bmoz.
\end{eqnarray}
Define the maximal operators $T^+$ and $\wz T^+$ as in Definition \ref{d2.4}
(i) with $S_t$ replaced by $T_t$ and $\wz T_t$, respectively.
  Then by \eqref{5.4}, there exists a positive constant $C$ such that
  for all $f\in\bmoz$,
$$\wz T^+(f)-C\|f\|_\bmoz\le T^+(f)\le \wz T^+(f)+C\|f\|_\bmoz.$$

 Since $\wz T^+$ is bounded from $\bmoz$ to $\bloz$
(see Theorem \ref{t5.1}), applying Lemma \ref{l4.2} again,
we have the following corollary.

\begin{cor}\hspace{-0.2cm}{\bf.}\hspace{0.2cm}\label{c5.2}
Assume that \eqref{5.3} holds. Then there exists a positive constant
$C$ such that for all
$f\in\bmo_{\rho}(\cx)$, $T^+(f)
\in\bmo_\rho(\cx)$ and $\|T^+(f)\|_{\blo_{\rho}(\cx)}
\le C\|f\|_{\bmo_{\rho}(\cx)}.$
\end{cor}

\noindent 5.2. {\bf Boundedness of Poisson semigroup maximal functions.}
Let $\{T_t\}_{t>0}$ be a family of bounded linear
integral operators on $\lt$ and
\begin{equation}\label{5.5}
P_t(f)\equiv\frac1{\sqrt\pi}\int_0^\fz \frac{e^{-s}}{\sqrt
s}T_{t/(2\sqrt s)}(f)\,ds.
\end{equation}
Define the maximal operator $P^+$ as in Definition
\ref{d2.4} (i) by replacing $T_t$ with $P_t$. If $\{T_t\}_{t>0}$ is replaced
by another family $\{\wz T_t\}_{t>0}$ of
bounded linear integral operators on $\lt$, we then denote the corresponding
$P_t$ in  \eqref{5.5} by $\wz P_t$ and the corresponding maximal operator by $\wz P^+$.

\begin{lem}\hspace{-0.2cm}{\bf.}\hspace{0.2cm}\label{l5.1}
(i) Let $\{\wz T_t\}_{t>0}$ be an $(\ez_1,\, \ez_2)$-$\ati$ for
$\ez_1\in(0, 1]$ and $\ez_2\in(0,\fz)$. Then $\{\wz P_t\}_{t>0}$ is
an $(\ez_1,\, \ez_2')$-$\ati$ with
$\ez_2'\in(0,\,\ez_2]\cap(0,\,1)$.

(ii) If $\{\wz T_t\}_{t>0}$ and $\{T_t\}_{t>0}$ satisfy \eqref{5.3}
with constants $C\in(0,\,\fz)$,
$\ez_1\in(0,\,1],\,\ez_2\in(0,\,\fz),\,\dz\in(0,\,1]$ and
$\gz\in(0,\,\fz)$, then so do $\{\wz P_t\}_{t>0}$ and
$\{P_t\}_{t>0}$ with constants $C'\in(0,\,\fz)$,
$\ez_1,\,\ez_2',\,\dz'$ and $\gz'$, where
$\ez_2'\in(0,\,\ez_2]\cap(0,\,1)$, $\gz'\in(0,\,\gz)$ and
$\dz'\in(0,\,\dz)$ satisfying $0<\gz'+\dz'<1$.
 \end{lem}

\begin{pf}\rm We first prove (i).
By $\int_\cx \wz T_t(x,\,y)\,d\mu(y)=1$ for all $x\in \cx$, we obtain
$$\int_\cx \wz P_t(x,\,y)\,d\mu(y)=
\frac1{\sqrt\pi}\int_\cx\int_0^\fz \frac{e^{-s}}{\sqrt s}\wz
T_{t/(2\sqrt s)}(x, y)\,ds \,d\mu(y)=\frac1{\sqrt\pi}\int_0^\fz
\frac{e^{-s}}{\sqrt s}\,ds=1.$$ Similarly, for any $x\in\cx$,
it follows from $\int_\cx
\wz T_t(y,\,x)\,d\mu(y)=1$  that $\int_\cx \wz
P_t(y,\,x)\,d\mu(y)=1$.

For all $s$, $t\in(0,\fz)$ and $x$, $y\in
\cx$, from the fact that
\begin{equation}\label{5.6}
t+d(x, y)\le (1+s)[t/s+d(x, y)],
\end{equation}
it follows that
\begin{equation}\label{5.7}
\frac{t/s}{t/s+d(x,\,y)}\le (1+s^{-1})\frac {t}{t+d(x,\,y)}.
\end{equation}
On the other hand, by \eqref{2.1} and \eqref{5.6},
we have that for all $s$, $t\in(0, \fz)$ and $x$, $y\in \cx$,
\begin{eqnarray}\label{5.8}
V_{t/s}(x)+V(x,\,y)&&\sim \mu(B(x,\,t/s+d(x,\,y)))\nonumber\\
&&\gs(1+s)^{- n } \mu(B(x,\,t+d(x,\,y)))\sim (1+s)^{- n
}[V_{t}(x)+V(x,\,y)].
\end{eqnarray}
Since $\{\wz T_t\}_{t>0}$ is an $(\ez_1, \ez_2)$-$\ati$,  by Definition \ref{d2.3} (i), \eqref{5.7} and
\eqref{5.8}, we obtain that for all $\ez'_2\in(0, \ez_2]$ and all $x,\,y\in\cx$,
\begin{eqnarray}\label{5.9}
\lf|\wz P_t(x,\,y)\r|&&\ls\int_0^\fz e^{-s^2/4}\lf|\wz T_{t/s}(x,\,y)\r|\,ds\nonumber\\
&&\ls \int_0^\fz  e^{-s^2/4} \frac1{V_{t/s}(x)+V(x,\,y)}\lf(\frac
{t/s}{t/s+d(x,\,y)}\r)^{\ez_2} \,ds\nonumber\\
&&\ls\frac1{V_{t}(x)+V(x,\,y)}\lf(\frac{t}{t+d(x,\,y)}\r)^{\ez_2'}.
\end{eqnarray}

Now we prove that for all $x$, $x'$, $y\in\cx$ and $t\in(0, \fz)$, if
$d(x,\,x')\le\frac12[t+d(x,\,y)]$, then
\begin{equation}\label{5.10}
\lf|\wz P_t(x,\,y)-\wz P_t(x',\,y)\r| \ls \lf(\frac
{d(x,\,x')}{t+d(x,\,y)}\r)^{\ez_1}\frac1{V_{t}(x)+V(x,\,y)}\lf(\frac
{t}{t+d(x,\,y)}\r)^{\ez_2'}.
\end{equation}
If $\frac18[t+d(x,\,y)]\le d(x,\,x')\le\frac12[t+d(x,\,y)]$,
\eqref{5.9} yields \eqref{5.10}.
If $\frac14d(x,\,y)\le d(x,\,x')<\frac18[t+d(x,\,y)]$, then
$d(x,\,y)\le 4 d(x,\,x')\le t$. In this case,
$t+d(x,\,y)\sim t+d(x',\,y)$.  By
Definition \ref{d2.3} (ii), \eqref{5.6}, \eqref{5.7} and
\eqref{5.8}, we obtain
\begin{eqnarray*}
&&\lf|\wz P_t(x,\,y)-\wz P_t(x',\,y)\r|\\
&&\hs\ls\int_0^\fz e^{-s^2/4} \lf|\wz T_{t/s}(x,\,y)-\wz T_{t/s}(x',\,y)\r|\,ds\\
&&\hs\ls \lf[\int_0^{t/[2d(x,\,x')]} \lf(\frac
{d(x,\,x')}{t/s+d(x,\,y)}\r)^{\ez_1}+ \int_{t/[2d(x,\,x')]}^\fz\r]e^{-s^2/4}
\frac1{V_{t/s}(x)+V(x,\,y)}\\
&&\hs\hs\times \lf(\frac{t/s}{t/s+d(x,\,y)}\r)^{\ez_2}\,ds\\
&&\hs\ls \lf(\frac
{d(x,\,x')}{t+d(x,\,y)}\r)^{\ez_1}\frac1{V_{t}(x)+V(x,\,y)}\lf(\frac
{t}{t+d(x,\,y)}\r)^{\ez_2'}\lf[\int_0^{t/[2d(x,\,x')]}(1+s)^{\ez_1}\r.\\
&&\hs\hs\lf.+ \int_{t/[2d(x,\,x')]}^\fz s^{\ez_1} \r]e^{-s^2/4}
(1+s)^ n (1+s^{-\ez_2'})\,ds\\
&&\hs\ls \lf(\frac
{d(x,\,x')}{t+d(x,\,y)}\r)^{\ez_1}\frac1{V_{t}(x)+V(x,\,y)}\lf(\frac
{t}{t+d(x,\,y)}\r)^{\ez_2'}.
\end{eqnarray*}

If $d(x,\,x')<\frac14d(x,\,y)$,
then by Definition \ref{d2.3} (ii), \eqref{5.6}, \eqref{5.7} and \eqref{5.8},
we have
\begin{eqnarray*}
&&|\wz P_t(x,\,y)-\wz P_t(x',\,y)|\\
&&\hs\ls\int_0^\fz e^{-s^2/4}\lf|\wz T_{t/s}(x,\,y)-\wz T_{t/s}(x',\,y)\r|\,ds\\
&&\hs\ls \int_0^\fz e^{-s^2/4}\lf(\frac
{d(x,\,x')}{t/s+d(x,\,y)}\r)^{\ez_1}\frac1{V_{t/s}(x)+V(x,\,y)}\lf(\frac
{t/s}{t/s+d(x,\,y)}\r)^{\ez_2}\,ds\\
&&\hs\ls \lf(\frac
{d(x,\,x')}{t+d(x,\,y)}\r)^{\ez_1}\frac1{V_{t}(x)+V(x,\,y)}\lf(\frac
{t}{t+d(x,\,y)}\r)^{\ez_2'},
\end{eqnarray*}which verifies \eqref{5.10}. Similarly with $x$ and $y$
interchanged, we have that $\wz P_t(x,\,y)$ satisfies Definition
\ref{d2.3} (iii). This shows (i).

To prove (ii), by (i), we only need to prove \eqref{5.3} for $\wz
P_t$ and $P_t$.
By \eqref{5.3}, \eqref{5.6},  \eqref{5.7} with $d(x,\,y)$ replaced
by $\rho(x)$, and \eqref{5.8}, we have that for all $\gz'\in(0,\,\gz)$ and
$\dz'\in(0,\,\dz')$ satisfying $0<\gz'+\dz'<1$ and $x$, $y\in\cx$,
\begin{eqnarray*}
&&|\wz P_t(x,\,y)-P_t(x,\,y)|\\
&&\hs\ls\int_0^\fz\frac{e^{-s}}{\sqrt s}
\lf|\wz T_{t/{(2\sqrt{s})}}(x,\,y)-T_{t/{(2\sqrt{s})}}(x,\,y)\r|\,ds\\
&&\hs\ls \int_0^\fz  e^{-4s^2} \lf(\frac
{t/s}{t/s+\rho(x)}\r)^{\dz}\frac1{V_{t/s}(x)+V(x,\,y)}\lf(\frac
{t/s}{t/s+d(x,\,y)}\r)^{\gz} \,ds\\
&&\hs\ls \lf(\frac{t}{t+\rho(x)}\r)^{\dz'}
\frac1{V_{t}(x)+V(x,\,y)}\lf(\frac{t}{t+d(x,\,y)}\r)^{\gz'}.
\end{eqnarray*}
This finishes the proof of Lemma \ref{l5.1}.
\hfill$\Box$\end{pf}

Then by Theorem \ref{t5.1}, Corollaries \ref{c5.1} and \ref{c5.2},
we obtain the following result.

\begin{thm}\hspace{-0.2cm}{\bf.}\hspace{0.2cm}\label{t5.2} Let $\rho$ be an
admissible function. Assume that $\{\wz T_t\}_{t>0}$ and
$\{T_t\}_{t>0}$ satisfy \eqref{5.3} with constants $C\in(0,\,\fz)$,
$\ez_1\in(0,1]$, $\ez_2\in(0,\,\fz)$, $\dz\in(0,\,1]$ and
$\gz\in(0,\,\fz)$. Then there exists a positive constant $C$ such
that for all $f\in\bmo_{\rho}(\cx)$, $P^+(f)$, $\wz P^+(f),\,\wz
P_\rho^+(f)\in\bloz$ and
$$\|P^+(f)\|_{\blo_{\rho}(\cx)}+\|\wz P^+(f)\|_{\blo_{\rho}(\cx)}+
\|\wz P_\rho^+(f)\|_{\blo_{\rho}(\cx)} \le
C\|f\|_{\bmo_{\rho}(\cx)}.$$
\end{thm}

\vspace{-1cm}
\section*{}
\setcounter{section}{6}
\setcounter{thm}{0}
\setcounter{defn}{0}
\setcounter{lem}{0}
\setcounter{rem}{0}
\setcounter{prop}{0}
\setcounter{cor}{0}
\noindent 6. {\bf Boundedness of the Littlewood-Paley
$g$-function.}  In this section, we consider the boundedness
of certain variant of the Littlewood-Paley $g$-function on $\bmoz$.

Let $\rho$ be an admissible function on $\cx$ and $\{Q_t\}_{t>0}$ be
a family of  operators bounded on $\lt$ with  integral kernels
$\{Q_t(x,\,y)\}_{t>0}$ satisfying that there exist constants
$C\in(0,\,\fz)$, $\dz_1\in(0,\,\fz)$,  $\bz\in(0,\,1]$,
$\dz_2\in(0,\,1)$ and $\gz\in(0,\,\fz)$ such that for all
$t\in(0,\,\fz)$ and $x,\,x',\,y\in\cx$ with $d(x,\,x')\le \frac t2$,

$(Q)_{\rm i}$ $|Q_t(x,\,y)|\le C\frac1{V_t(x)+V(x,\,y)}(\frac
t{t+d(x,\,y)})^{\gz}(\frac {\rho(x)}{t+\rho(x)})^{\dz_1}$;

$(Q)_{\rm ii}$  $|Q_t(x,\,y)-Q_t(x',\,y)|\le C(\frac
{d(x,\,x')}{t+d(x,\,y)})^\bz \frac{1}{V_t(x)+V(x,\,y)}(\frac
t{t+d(x,\,y)})^{\gz}$;

$(Q)_{\rm iii}$ $|\int_\cx Q_t(x,\,y)d\mu(y)|\le C(\frac
t{t+\rho(x)})^{\dz_2}$.

For all $f\in L^1_\loc(\cx)$ and $x\in\cx$, define the Littlewood-Paley
$g$-function by
\begin{equation}\label{6.1}
g(f)(x)\equiv\lf(\int_0^\fz |Q_t(f)(x)|^2\frac {dt}t\r)^{1/2}.
\end{equation}

\begin{thm}\hspace{-0.2cm}{\bf.}\hspace{0.2cm}\label{t6.1}
Let $\rho$ be an admissible function on $\cx$. Suppose the
$g$-function defined in \eqref{6.1} is bounded on $\lt$. Then there
exists a positive constant $C$ such that for all
$f\in\bmo_{\rho}(\cx)$, $[g(f)]^2\in\blo_{\rho}(\cx)$ and
$\|[g(f)]^2\|_{\blo_{\rho}(\cx)}\le C\|f\|^2_{\bmo_{\rho}(\cx)}$.
\end{thm}

\begin{pf}\rm
By the homogeneity of $\|\cdot\|_\bmoz$ and $\|\cdot\|_\bloz$, we
assume that $f\in\bmo_{\rho}(\cx)$ and $\|f\|_{\bmo_{\rho}(\cx)}=1.$
We first prove that for all balls $B\equiv B(x_0,\,r)$ with
$r\ge\rho(x_0)$,
\begin{equation}\label{6.2}
\frac1{\mu(B)}\int_B[g(f)(x)]^2\,d\mu(x)\ls1.
\end{equation}
For any $x\in B$, write
$$[g(f)(x)]^2\equiv\int_0^{8\rho(x)} |Q_t(f)(x)|^2\frac {dt}t
+\int_{8\rho(x)}^\fz |Q_t(f)(x)|^2\frac {dt}t\equiv
[g_1(f)(x)]^2+[g_2(f)(x)]^2.$$

By the $L^2(\cx)$-boundedness of $g$,
\eqref{2.1} and Lemma \ref{l3.1}, we have
\begin{equation}\label{6.3}
\frac1{\mu(B)}\int_B[g_1(f\chi_{2B})(x)]^2\,d\mu(x)\ls
\frac1{\mu(B)}\int_{2B}|f(x)|^2\,d\mu(x)\ls1.
\end{equation}
For any $x\in B$, by ${\rm (Q)_i}$,
\begin{eqnarray*}
\lf|Q_t\lf(f\chi_{(2B)^\complement}\r)(x)\r|&&\ls\int_{(2B)^\complement}
\frac1{V_t(x)+V(x,\,y)}\lf(\frac{t}{t+d(x,\,y)}\r)^{\gz}
|f(y)|\,d\mu(y)\\
&&\ls\lf(\frac tr\r)^{\gz}\sum_{j=1}^\fz
\frac{2^{-j\gz}}{V_{2^{j-1}r}(x)}
\int_{d(x,\,y)<2^jr}|f(y)|\,d\mu(y) \ls\lf(\frac
tr\r)^{\gz}.\nonumber
\end{eqnarray*}
Notice that for all $x\in B$, by \eqref{3.7}, we have $\rho(x)\ls r$.
From the inequality above we deduce that
\begin{equation}\label{6.4}
\frac1{\mu(B)}\int_B\lf[g_1
\lf(f\chi_{(2B)^\complement}\r)(x)\r]^2\,d\mu(x)\ls
\frac1{\mu(B)}\int_B\int_0^{8\rho(x)}\lf(\frac tr\r)^{2\gz}
\,\frac{dt}t\,d\mu(x)\ls1.
\end{equation}
Combining \eqref{6.3} and \eqref{6.4} gives us that
\begin{equation}\label{6.5}
\frac1{\mu(B)}\int_B[g_1(f)(x)]^2\,d\mu(x)\ls1.
\end{equation}

To prove \eqref{6.2} with $g_2$,  we first
 notice that for all $x\in B$ and $t\ge8\rho(x)$,
 \begin{eqnarray}\label{6.6}
|Q_t(f)(x)|&&\ls\int_{\cx}\frac1{V_t(x)+V(x,\,y)}
\lf(\frac{t}{t+d(x,\,y)}\r)^{\gz}\lf(\frac{\rho(x)}{t+\rho(x)}\r)^{\dz_1}
|f(y)|\,d\mu(y)\nonumber\\
&&\ls\lf(\frac{\rho(x)}{t}\r)^{\dz_1}\sum_{j=0}^\fz
2^{-j\gz}\frac1{V_{2^{j-1}t}(x)}
\int_{d(x,\,y)<2^jt}|f(y)|\,d\mu(y)\ls\lf(\frac{\rho(x)}{t}\r)^{\dz_1}.
\end{eqnarray}
Then
\begin{equation*}
\frac1{\mu(B)}\int_B[g_2(f)(x)]^2\,d\mu(x)\ls
\frac1{\mu(B)}\int_B\int_{8\rho(x)}^\fz
\lf(\frac{\rho(x)}t\r)^{2\dz_1}\frac {dt}t\,d\mu(x)\ls1,
\end{equation*}
which together with \eqref{6.5} gives \eqref{6.2}. Moreover, since
\eqref{6.2} holds for all balls $B(x_0,\,r)$ with $r\ge\rho(x_0)$,
we have that $g(f)(x)<\fz$ for a.\,e. $x\in\cx$.

Now we assume that $B\equiv B(x_0,\,r)$ with $r<\rho(x_0)$. If $r\ge
\rho(x_0)/8$, then by \eqref{2.1} and \eqref{6.2}, we have
\begin{equation*}
\frac1{\mu(B)}\int_B\lf\{[g(f)(x)]^2-{\mathop\einf_B}[g(f)]^2\r\}\,d\mu(x)
\ls \frac1{\mu(8B)}\dint_{8B}[g(f)(x)]^2\,d\mu(x)\ls1,
\end{equation*}
which is desired. Assume that $r<\rho(x_0)/8$.
It suffices to show that for $\mu$-a.\,e.\,$y\in B$,
\begin{eqnarray*}
&&\frac1{\mu(B)}\int_B\lf\{[g(f)(x)]^2-[g(f)(y)]^2\r\}\,d\mu(x)\ls1.
\end{eqnarray*}
For all $x\in B$, write
\begin{eqnarray*}
[g(f)(x)]^2&&=\int_0^{8r} |Q_t(f)(x)|^2\frac {dt}t
+\int_{8r}^{8\rho(x_0)}\cdots+\int_{8\rho(x_0)}^\fz\cdots\\
&&\equiv[g_r(f)(x)]^2+[g_{r,\,x_0}(f)(x)]^2+[g_\fz(f)(x)]^2.
\end{eqnarray*}
Observe that for $\mu$-a.\,e.\,$y\in B$,
\begin{eqnarray*}
&&\frac1{\mu(B)}\int_B\lf\{[g(f)(x)]^2-[g(f)(y)]^2\r\}\,d\mu(x)\\
&&\hs\le \frac1{\mu(B)}\int_B\lf\{[g_r(f)(x)]^2+[g_\fz(f)(x)]^2+
[g_{r,\,x_0}(f)(x)]^2-[g_{r,\,x_0}(f)(y)]^2 \r\}\,d\mu(x).
\end{eqnarray*}

We first prove that
\begin{equation}\label{6.7}
\frac1{\mu(B)}\int_B[g_r(f)(x)]^2\,d\mu(x)\ls1.
\end{equation}
Write $f\equiv f_1+f_2+f_B,$ where $f_1\equiv(f-f_B)\chi_{2B}$ and
$f_2\equiv(f-f_B)\chi_{(2B)^\complement}$. By the
$L^2(\cx)$-boundedness of $g$, \eqref{2.1} and Lemma \ref{l3.1}, we
have
\begin{equation}\label{6.8}
\frac1{\mu(B)}\int_B[g_r(f_1)(x)]^2\,d\mu(x) \ls
\frac1{\mu(B)}\int_{2B}|f-f_B|^2\,d\mu(x)\ls1.
\end{equation}
For all $x\in B$, by ${\rm (Q)_i}$, \eqref{2.1} and the
fact that $|f_{2^{j+1}B}-f_B|\ls j$ for all $j\in\nn$,
we have
\begin{eqnarray*}
|Q_t(f_2)(x)|&&\ls\int_{(2B)^\complement}
\frac1{V_t(x)+V(x,\,z)}\lf(\frac{t}{t+d(x,\,z)}\r)^{\gz}
|f(z)-f_B|\,d\mu(z)\\
&&\ls\sum_{j=1}^\fz \lf(\frac
t{2^{j-1}r}\r)^{\gz}\lf[\frac1{V_{2^{j-1}r}(x)}
\int_{2^{j+1}B}\lf[|f(z)-f_{2^{j+1}B}|+|f_{2^{j+1}B}-f_B|\r]
\,d\mu(z)\r]\nonumber\\
&&\ls\lf(\frac tr\r)^{\gz}\sum_{j=1}^\fz j2^{-j\gz} \ls\lf(\frac
tr\r)^{\gz},\nonumber
\end{eqnarray*}
which further implies that
\begin{equation}\label{6.9}
\frac1{\mu(B)}\int_B[g_r(f_2)(x)]^2\,d\mu(x)\ls \int_0^{8r}\lf(\frac
tr\r)^{2\gz} \,\frac{dt}t\ls1.
\end{equation}
By \eqref{6.8} and \eqref{6.9}, to prove \eqref{6.7}, it remains to
show that
\begin{equation}\label{6.10}
\frac1{\mu(B)}\int_B[g_r(f_B)(x)]^2\,d\mu(x)\ls1.
\end{equation}
Observe that $|f_B|\ls \log\frac{\rho(x_0)}{r}$.
 For all $x\in B$ and $t\in (0,\rho(x_0))$, from
  ${\rm (Q)_{iii}}$ and the fact that $r<\rho(x_0)/8$ and \eqref{3.7},
  it follows that
$$|Q_t(f_B)(x)|\ls \lf(\frac t{\rho(x)}\r)^{\dz_2}|f_B|\ls \lf(\frac
t{\rho(x_0)}\r)^{\dz_2} \log\frac{\rho(x_0)}{r},$$ which via $t\le
8r<\rho(x_0)$ further yields \eqref{6.10}.

Let $a\in[1/8, \fz)$ and ${\wz C}_a$ be as in \eqref{3.7}.
We now claim that for all
$f\in \bmoz$ with $\|f\|_\bmoz=1$, $x\in B$ and $t\le
8{\wz C}_a\rho(x_0)$,
\begin{equation}\label{6.11}
|Q_t(f)(x)|\ls 1.
\end{equation}
 In fact, by ${\rm (Q)_{ii}}$, we obtain
\begin{eqnarray}\label{6.12}
&&|Q_t(f-f_{B(x,\,t)})(x)|\nonumber\\
&&\hs\ls\int_{\cx}\frac1{V_t(x)+V(x,\,y)}
\lf(\frac{t}{t+d(x,\,y)}\r)^{\gz}
|f(y)-f_{B(x,\,t)}|\,d\mu(y)\nonumber\\
&&\hs\ls\sum_{j=0}^\fz 2^{-j\gz}\frac1{V_{2^{j-1}t}(x)}
\int_{d(x,\,y)<2^jt}|f(y)-f_{B(x,\,t)}|\,d\mu(y)\ls1.
\end{eqnarray}
Since $\dz_2>0$, by  ${\rm (Q)_{iii}}$ and the fact that for all $x\in
\cx$, $|f_{B(x,\,t)}|\ls 1+\log\frac{\rho(x)}{t}$, we have
\begin{eqnarray*}
|Q_t(f_{B(x,\,t)})(x)|\ls\lf(\frac t{\rho(x)}\r)^{\dz_2}
\lf(1+\log\frac{\rho(x)}{t}\r)\ls1.
\end{eqnarray*}
Combining this and \eqref{6.12} proves the claim.

Using \eqref{6.11}, \eqref{3.7} and
\eqref{6.6}, we have that for all $x\in B$,
\begin{eqnarray*}
\int_{8\rho(x_0)}^\fz|Q_t(f)(x)|^2\,\frac{dt}t
&\le&\int_{8{\wz C}_a\rho(x_0)}^\fz|Q_t(f)(x)|^2\,\frac{dt}t
+\int_{8\rho(x_0)}^{8{\wz C}_a\rho(x_0)}\cdots\\
&\ls&\int_{8{\wz C}_a\rho(x_0)}^\fz\lf(\frac{\rho(x)}{t}\r)^{2\dz_1}
\,\frac{dt}t+1\ls1,
\end{eqnarray*}
which gives that
\begin{equation}\label{6.13}
\frac1{\mu(B)}\int_B[g_\fz(f)(x)]^2\,d\mu(x)\ls1.
\end{equation}

By \eqref{6.7} and \eqref{6.13}, to finish the proof of Theorem \ref{t6.1},
it remains to show that for
$\mu$-a.\,e. $y\in B$,
\begin{equation}\label{6.14}
\frac1{\mu(B)}\int_B\lf\{[g_{r,\,x_0}(f)(x)]^2-[g_{r,\,x_0}(f)(y)]^2
\r\}\,d\mu(x)\ls1.
\end{equation}

From \eqref{6.11}, we deduce that for $\mu$-a.\,e.\, $x$, $y\in B$,
\begin{eqnarray*}
&&\lf\{[g_{r,\,x_0}(f)(x)]^2-[g_{r,\,x_0}(f)(y)]^2\r\}\\
&&\hs\le\int_{8r}^{8\rho(x_0)}
|Q_t(f)(x)+Q_t(f)(y)||Q_t(f)(x)-Q_t(f)(y)|\,\frac{dt}t\\
&&\hs\ls\int_{8r}^{8\rho(x_0)}|Q_t(f)(x)-Q_t(f)(y)|\,\frac{dt}t.
\end{eqnarray*}
For $t\in(8r, 8\rho(x_0))$ and $x$, $y\in B$, we write
\begin{eqnarray*}
&&|Q_t(f)(x)-Q_t(f)(y)|\\
&&\hs\le \lf|\dint_\cx\lf[Q_t(x, z)-Q_t(y,
z)\r][f(z)-f_B]\,d\mu(z)\r|
+|f_B|\lf|\dint_\cx\lf[Q_t(x, z)-Q_t(y, z)\r]\,d\mu(z)\r|\\
&&\hs\equiv{\rm J_1}+{\rm J_2}.
\end{eqnarray*}
By  ${\rm (Q)_{ii}}$, $t\in (8r, 8\rho(x_0))$, \eqref{2.1} and the
fact that $2^{j+1}B\subset B(x,\,2^{j+2}r)$ for all $x\in B$,
we obtain
\begin{eqnarray*}
{\rm J_1}&&\ls\int_\cx
\lf(\frac{d(x,\,y)}{t+d(x,\,z)}\r)^\bz\frac1{V_t(x)+V(x,\,z)}
\lf(\frac{t}{t+d(x,\,z)}\r)^{\gz}|f(z)-f_B|\,d\mu(z)\\
&&\ls\sum_{j=0}^\fz\frac{r^\bz t^{\gz}}{(t+2^{j-1}r)^{\bz+\gz}}
\lf\{\frac1{V_{2^{j-1}r}(x)}\int_{2^{j+1}B}|f(z)-f_{2^{j+1}B}|
\,d\mu(z)+|f_{2^{j+1}B}-f_B|\r\}\\
&&\ls\sum_{j=0}^\fz\frac{r^\bz t^{\gz}}{(t+2^jr)^{\bz+\gz}}(j+2).
\end{eqnarray*}
From this, we deduce that
\begin{equation}\label{6.15}
\dint_{8r}^{8\rho(x_0)}{\rm J_1}\,\frac{dt}t \ls
r^\bz\sum_{j=0}^\fz(j+2)\dint_{8r}^{8\rho(x_0)}
\frac{t^{\gz-1}}{(t+2^jr)^{\bz+\gz}}\,dt\ls 1.
\end{equation}

On the other hand, we obtain that for $\mu$-a.\,e.\, $x$, $y\in B$,
\begin{eqnarray*}
&&\lf|\dint_\cx\lf[Q_t(x, z)-Q_t(y, z)\r]\,d\mu(z)\r|\\
&&\hs\ls\int_\cx
\lf(\frac{d(x,\,y)}{t+d(x,\,z)}\r)^\bz\frac1{V_t(x)+V(x,\,z)}
\lf(\frac{t}{t+d(x,\,z)}\r)^{\gz}\,d\mu(z)\\
&&\hs\ls\sum_{j=0}^\fz\lf(\frac r{2^jt}\r)^\bz2^{-j\gz}\ls
\lf(\frac r{t}\r)^\bz.
\end{eqnarray*}
The fact that $r<\rho(x_0)/8$ together with \eqref{3.7} and
$|f_B|\ls \log\frac{\rho(x_0)}r$ implies that for $\mu$-a.\,e.\,$x$,
$y\in B$,
\begin{eqnarray*}
\dint_{8r}^{8\rho(x_0)}{\rm J_2}\,\frac{dt}t
&\le&\dint_{8r}^{8\rho(x_0)}\log\frac{\rho(x_0)}r
\lf|\dint_\cx\lf[Q_t(x, z)-Q_t(y, z)\r]\,d\mu(z)\r|^{\frac23}
\lf(\frac{t}{\rho(x_0)}\r)^{\frac{\dz_2}3}\,\frac{dt}t\\
&\ls &\dint_{8r}^{8\rho(x_0)}\log\frac{\rho(x_0)}r \lf(\frac
r{t}\r)^{\frac{\bz}3}
\lf(\frac{r}{\rho(x_0)}\r)^{\min\lf(\frac{\bz}3,
\frac{\dz_2}3\r)}\,\frac{dt}t
\ls\dint_{8r}^{\fz}\lf(\frac r{t}\r)^{\frac{\bz}3}\,\frac{dt}t
\ls1.
\end{eqnarray*}
This together with \eqref{6.15} leads to \eqref{6.14}, and hence,
finishes the proof of Theorem \ref{t6.1}.
\hfill$\Box$\end{pf}

As a consequence of Theorem \ref{t6.1}, we have the following
conclusion.

\begin{cor}\hspace{-0.2cm}{\bf.}\hspace{0.2cm}\label{c6.1}
With the assumptions same as in Theorem \ref{t6.1},
then there exists a positive constant $C$ such
that for all $f\in\bmo_{\rho}(\cx)$, $g(f)\in\blo_{\rho}(\cx)$ and
$\|g(f)\|_{\blo_{\rho}(\cx)} \le C\|f\|_{\bmo_{\rho}(\cx)}.$
\end{cor}

\begin{pf}\rm
Since
\begin{equation*}
g(f)-{\mathop\einf_Bg(f)}\le\lf\{[g(f)]^2-
{\mathop\einf_B[g(f)]^2}\r\}^{1/2},
\end{equation*}
by the H\"older inequality and Theorem \ref{t6.1},
we have that for all balls $B\notin\cd$,
\begin{eqnarray*}
&&\frac1{\mu(B)}\int_B \lf[g(f)(x)-{\mathop\einf_{B}}g(f)\r]\,d\mu(x)\\
&&\hs\ls \frac1{\mu(B)}\dint_B \lf\{\lf[g(f)(x)\r]^2-
{\mathop\einf_B}[g(f)]^2\r\}^{1/2}\,d\mu(x)\\
&&\hs\ls \lf\{\frac1{\mu(B)}\dint_B \lf\{[g(f)(x)]^2-
{\mathop\einf_B}[g(f)]^2\r\}\,d\mu(x)\r\}^{1/2}\ls\|f\|_{\bmo_{\rho}(\cx)}.
\end{eqnarray*}
On the other hand, by \eqref{6.2} and the H\"older inequality,
we obtain that for all balls $B\in\cd$,
\begin{equation*}
\frac1{\mu(B)}\int_Bg(f)(x)\,d\mu(x)\le
\lf(\frac1{\mu(B)}\int_B[g(f)(x)]^2\,d\mu(x)\r)^{1/2}\ls\|f\|_{\bmo_{\rho}(\cx)}.
\end{equation*}
Combining the two inequalities above finishes the proof of Corollary
\ref{c6.1}.
\hfill$\Box$\end{pf}

\vspace{-1.4cm}
\section*{}
\setcounter{section}{7}
\setcounter{thm}{0}
\setcounter{defn}{0}
\setcounter{lem}{0}
\setcounter{rem}{0}
\setcounter{prop}{0}
\setcounter{cor}{0}
\noindent 7. {\bf Applications.}
This section is divided into Subsections 7.1 through 7.4,
which are devoted to the applications
of results obtained in Sections 5
and 6, respectively,
to the Schr\"odinger operator or the
degenerate Schr\"odinger operator on $\rd$, the
sub-Laplace Schr\"odinger operator on Heisenberg groups or on
connected and simply connected nilpotent Lie groups.

\medskip
\noindent 7.1. {\bf Schr\"odinger operators on $\rd$.}\label{s7.1}
Let $d$ be a positive integer and $d\ge3$, and $\rd$ be the
$d$-dimensional Euclidean space endowed with the Euclidean norm
$|\cdot|$ and the Lebesgue measure $dx$. Denote the Laplacian
$\sum_{j=1}^d\frac{\partial^2}{\partial x_j^2}$ on $\rd$ by
$\Delta$ and the corresponding heat (Gauss) semigroup
$\{e^{t\Delta}\}_{t>0}$ by $\{\wz T_t\}_{t>0}$. Let $V$ be a
nonnegative locally integrable function on $\rd$,
$\cl\equiv-\Delta+V$ be the Schr\"odinger operator and
$\{T_t\}_{t>0}$ be the corresponding semigroup.
Moreover, for all $t>0$ and $x,\,y\in\rd$, set
$$Q_t(x,\,y)\equiv t^2\frac{d T_s(x,\,y)}{ds}\Bigg|_{s=t^2}.$$
Let $q \in(d/2, d]$,
$V\in\cb_q(\rd, |\cdot|, dx)$ and $\rho$ be as in \eqref{2.4}.
Then we have the following estimates; see \cite{d05,dz02,dz03}.

\begin{prop}\hspace{-0.2cm}{\bf.}\hspace{0.2cm}\label{p7.1}
Let $q \in(d/2, d]$, $\bz\in(0, 2-d/q)$ and $N\in\nn$.
Then there exist  positive constants $\wz C$ and $C$,
where $C$ is independent of $N$, such that
for all $t\in(0,\,\fz)$ and $x,\,x',\,y\in\cx$ with $d(x,\,x')\le
\sqrt t/2$,
\begin{enumerate}
\vspace{-0.2cm}
\item[(i)] $|T_t(x,\,y)|\le \wz Ct^{-d/2}\exp\{-\frac
{|x-y|^2}{Ct}\}[\frac {\rho(x)}{\sqrt t+\rho(x)}]^N
[\frac {\rho(y)}{\sqrt t+\rho(y)}]^N$;
\vspace{-0.2cm}
\item[(ii)]  $|T_t(x,\,y)-T_t(x',\,y)|\le\wz C
 [\frac {|x-x'|}{\sqrt t}]^{\bz}t^{-d/2}\exp \{-\frac
{|x-y|^2}{Ct}\}[\frac {\rho(x)}{\sqrt t+\rho(x)}]^N
[\frac {\rho(y)}{\sqrt t+\rho(y)}]^N$;
\vspace{-0.2cm}
\item[(iii)]  $|T_t(x,\,y)-\wz T_t(x,\,y)|\le\wz C
 [\frac {\sqrt t}{\sqrt t+\rho(x)}]^{2-d/q}t^{-d/2}\exp\{-\frac
{|x-y|^2}{Ct}\}$,
\end{enumerate}
\vspace{-0.1cm}
\noindent and for all $t\in(0,\,\fz)$ and $x,\,x',\,y\in\cx$ with $d(x,\,x')\le
t/2$,
\begin{enumerate}
\vspace{-0.3cm}
\item[(iv)] $|Q_t(x,\,y)|\le  \wz C t^{-d}\exp \{-\frac
{|x-y|^2}{Ct^2}\}[\frac {\rho(x)}{t+\rho(x)}]^N[\frac {\rho(y)}{t+\rho(y)}]^N$;
\vspace{-0.3cm}
\item[(v)]  $|Q_t(x,\,y)-Q_t(x',\,y)|\le\wz C
 [\frac {|x-x'|}{t}]^{\bz}t^{-d}\exp \{-\frac
{|x-y|^2}{Ct^2}\}[\frac {\rho(x)}{t+\rho(x)}]^N
[\frac {\rho(y)}{t+\rho(y)}]^N$;
\vspace{-0.3cm}
\item[(vi)] $|\int_\rd Q_t(x,\,y)d\mu(y)|\le \wz C
[\frac t{\rho(x)}]^{2-d/q}[\frac{\rho(x)}{t+\rho(x)}]^N$.
\end{enumerate}
\end{prop}

Observe that $\{\wz T_{t^2}\}_{t>0}$ is a continuous $(1,
N)$-$\ati$ for all positive constants $N$. Thus $\{T_{t^2}\}_{t>0}$ and
$\{\wz T_{t^2}\}_{t>0}$ satisfy the assumption \eqref{5.3}.
Moreover, the $L^2(\rd)$-boundedness of $g$-function was obtained in \cite{d05}.
Using these facts and Proposition \ref{p7.1} and applying
Theorems \ref{t5.1}, \ref{t5.2} and \ref{t6.1},
and Corollaries \ref{c5.1}, \ref{c5.2} and \ref{c6.1},
we have the following result.

\begin{prop}\hspace{-0.2cm}{\bf.}\hspace{0.2cm}\label{p7.2}
Let $q\in(d/2,\fz]$,  $V\in \cb_q(\rd,\,|\cdot|,\,dx)$ and
  $\rho$ be as in \eqref{2.4}.
There exists a positive constant $C$ such that
for all $f\in\bmo_{\rho}(\rd)$,
 $T^+(f)$, $\wz T^+(f)$, $\wz T^+_\rho(f)$, $P^+(f)$, $\wz P^+(f)$,
 $\wz P^+_\rho(f),
 \,g(f),\,[g(f)]^2\in\blo_\rho(\rd)$
and
\begin{eqnarray*}
&&\|T^+(f)\|_{\blo_\rho(\rd)}
+\|\wz T^+(f)\|_{\blo_\rho(\rd)}+
\|T^+_\rho(f)\|_{\blo_\rho(\rd)}+\|P^+(f)\|_{\blo_\rho(\rd)}\\
&&\hs\hs+\|\wz P^+(f)\|_{\blo_\rho(\rd)}
+\|\wz P^+_\rho(f)\|_{\blo_\rho(\rd)}
+\|g(f)\|_{\blo_\rho(\rd)}
+\|[g(f)]^2\|_{\blo_\rho(\rd)}^{1/2}\\
&&\hs\le C\|f\|_{\bmo_\rho(\rd)}.
\end{eqnarray*}
\end{prop}

We also point out that when $\rho$ is as in \eqref{2.4},
Dziuba\'nski et al \cite{dgmtz05} obtained the boundedness of
$T^+$, $P^+$ and $g$ on $\bmo_\rho(\rd)$.
Proposition \ref{p7.2} improves their results.

\medskip
\noindent 7.2. {\bf Degenerate Schr\"odinger operators on  $\rd$.}
Let $d\ge 3$ and $\rd$ be the $d$-dimensional Euclidean space
endowed with the Euclidean norm $|\cdot|$ and the Lebesgue measure $dx$.
Recall that a nonnegative locally integrable function $w$ is said to
be an $A_2(\rd)$ weight in the sense of Muckenhoupt if
$$\sup_{B\subset \rd}\lf\{\frac1{|B|}\int_B w(x)\,dx\r\}^{1/2}
\lf\{\frac1{|B|}\int_B [w(x)]^{-1}\,dx\r\}^{1/2}<\fz,$$
where the supremum is taken over all the balls in $\rd$.
Observe that if we set $w(E)\equiv\int_Ew(x)dx$ for any measurable set $E$, then
 there exist positive constants $C,\,Q$ and $\kz$ such that for all $x\in\rd$, $\lz>1$ and $r>0$,
$$C^{-1}\lz^\kz w(B(x,\,r))\le w(B(x,\,\lz r))\le C\lz^Qw(B(x,\,r)),$$
namely, the measure $w(x)\,dx$ satisfies \eqref{2.1}.
Thus $(\rd,\,|\cdot|,\,w(x)\,dx)$ is an RD-space.

Let $w\in A_2(\rd)$ and
$\{a_{i,\,j}\}_{1\le i,\,j\le d}$ be a real symmetric matrix function satisfying that
for all $x,\,\xi\in\rd$,
$$C^{-1}|\xi|^{2}\le\sum_{1\le i,\,j\le d}a_{i,\,j}(x)\xi_i\overline \xi_j
\le C|\xi|^2.$$ Then the degenerate elliptic operator
$\cl_0$ is defined by
$$\cl_0 f(x)\equiv-\frac1{w(x)}\sum_{1\le i,\,j\le d}
\partial_i(a_{i,\,j}(\cdot)\partial_j f)(x),$$
where $x\in\rd$.
Denote by $\{\wz T_t\}_{t>0}\equiv \{e^{-t\cl_0}\}_{t>0}$
the semigroup generated by $\cl_0$. We also denote the kernel of
$\wz T_t$ by $\wz T_t(x,\,y)$ for all $x,\,y\in\rd$ and $t\in(0, \fz)$.
Then it is known that there exist positive constants $C,\,C_6,\,\wz C_6$
and $\az\in(0,\,1]$ such that
for all $t\in(0, \fz)$ and $x,\,y\in\rd$,
\begin{equation*}
C^{-1}\frac1{V_{\sqrt t}(x)}\exp\lf\{-\frac{|x-y|^2}{\wz C_6t}\r\}\le
\wz T_t(x,\,y)\le
C\frac1{V_{\sqrt t}(x)}\exp\lf\{-\frac{|x-y|^2}{C_6t}\r\};
\end{equation*}
that for all $t\in(0, \fz)$ and $x,\,y,\,y'\in\rd$ with $|y-y'|<|x-y|/4$,
\begin{equation*}
 |\wz T_t(x,\,y)-\wz T_t(x,\,y')|\le
C\frac1{V_{\sqrt t}(x)}\lf(\frac{|y-y'|}{\sqrt t}\r)^\az\exp
\lf\{-\frac{|x-y|^2}{C_6t}\r\};
\end{equation*}
and that for all $t\in(0, \fz)$ and $x,\,y\in\rd$,
\begin{equation*}
\int_\rd \wz T_t(x,\,z)\,w(z)\,dz=1=\int_\rd \wz T_t(z,\,y)\,w(z)\,dz;
\end{equation*}
see, for example, Theorems 2.1, 2.7, 2.3, 2.4 and Corollary 3.4 of \cite{hs01}.

Let $V$ be a nonnegative locally integrable function on $w(x)\,dx$.
Define the  degenerate Schr\"odinger operator by $\cl \equiv\cl_0+V.$
Then $\cl$ generates a semigroup $\{T_t\}_{t>0}\equiv\{e^{-t\cl}\}_{t>0}$ with
kernels $\{T_t(x,\,y)\}_{t>0}$.
Moreover, for all $t\in(0, \fz)$ and $x,\,y\in\rd$, set
$$Q_t(x,\,y)\equiv t^2\frac{d T_s(x,\,y)}{ds}\Bigg|_{s=t^2}.$$
Let $q \in(Q/2, Q]$,
$V\in\cb_q(\rd, |\cdot|, w(x)\,dx)$ and $\rho$ be as in \eqref{2.4}.
Then $\{T_t\}_{t>0}$ and $\{Q_t\}_{t>0}$ satisfy Proposition \ref{p7.1}
with $t^{-d/2}$ replaced by $V_{\sqrt t}(x)$,
$t^{-d}$ by $V_{ t}(x)$, and $d $ by $Q $.

In fact, the corresponding Proposition \ref{p7.1}
(i) and (iii) here were given in \cite{d05}.
The proof of (ii) here is similar to that of Proposition \ref{l7.1};
see \cite{dz03} and also Lemma \ref{l7.4} below.
The proofs of the corresponding Proposition \ref{p7.1}
(iv), (v) and (vi) here are similar to that of Proposition 4 of \cite{dgmtz05}.
We omit the details here.

Observe that $\{\wz T_{t^2}\}_{t>0}$ is a continuous $(1,
N)$-$\ati$ for all positive constants $N$. Thus $\{T_{t^2}\}_{t>0}$ and
$\{\wz T_{t^2}\}_{t>0}$ satisfy the assumption \eqref{5.3}.
Moreover, the $L^2(\rd)$-boundedness of
$g$-function can be obtained by the same argument
as in Lemma 3 of \cite{d05}.
Using these facts and applying Theorems \ref{t5.1}, \ref{t5.2} and \ref{t6.1},
and Corollaries \ref{c5.1}, \ref{c5.2} and \ref{c6.1},
we have the following conclusions.

\begin{prop}\hspace{-0.2cm}{\bf.}\hspace{0.2cm}\label{p7.3}
Let $w\in A_2(\rd)$. Let
$q\in (Q/2,\fz]$,  $V\in \cb_q(\rd,\,|\cdot|,\,w(x)\,dx)$
and $\rho$
be as in \eqref{2.4} with $d\mu=w(x)\,dx$.
Then there exists a positive constant $C$ such that
for all $f\in\bmo_{\rho}(w(x)\,dx)$,
 $T^+(f)$, $\wz T^+(f)$, $\wz T^+_\rho(f)$, $P^+(f)$, $\wz P^+(f)$,
 $\wz P^+_\rho(f)$, $g(f)$, $[g(f)]^2\in\blo_{\rho}(w(x)\,dx)$ and
\begin{eqnarray*}
&&\|T^+(f)\|_{\blo_{\rho}(w(x)\,dx)}
+\|\wz T^+(f)\|_{\blo_{\rho}(w(x)\,dx)}+
\|T^+_\rho(f)\|_{\blo_{\rho}(w(x)\,dx)}\\
&&\hs\hs+\|P^+(f)\|_{\blo_{\rho}(w(x)\,dx)}+\|\wz P^+(f)\|_{\blo_{\rho}(w(x)\,dx)}
+\|\wz P^+_\rho(f)\|_{\blo_{\rho}(w(x)\,dx)}\\
&&\hs\hs+\|g(f)\|_{\blo_{\rho}(w(x)\,dx)}+
\|[g(f)]^2\|_{\blo_{\rho}(w(x)\,dx)}^{1/2}\\
&&\hs\le C\|f\|_{\bmo_{\rho}(w(x)\,dx)}.
\end{eqnarray*}
\end{prop}

\medskip
\noindent 7.3. {\bf Schr\"odinger operators on Heisenberg groups.}
The $(2n+1)$-dimensional Heisenberg group $\hh^n$ is a
connected and simply connected nilpotent Lie groups
with the underlying manifold $\rr^{2n}\times \rr$ and the multiplication
$$(x,\,s)(y,\,s)=\lf(x+y,\,t+s+2\sum_{j=1}^n[x_{n+j}y_j-x_jy_{n+j}]\r).$$
The homogeneous norm on $\hh^n$ is defined by
$|(x,\,t)|=(|x|^4+|t|^2)^{1/4}$ for all $(x,\,t)\in\hh^n$, which
induces a left-invariant metric
$d((x,\,t),\,(y,\,s))=|(-x,\,-t)(y,\,s)|$. Moreover, there exists a
positive constant $C$ such that $|B((x,\,t),\,r)|=Cr^Q,$ where
$Q=2n+2$ is the homogeneous dimension of $\hh^n$ and
$|B((x,\,t),\,r)|$ is the Lebesgue measure of the ball
$B((x,\,t),\,r)$. The triplet $(\hh^n,\, d,\,dx)$ is an RD-space.

A basis for the Lie algebra of left invariant vector fields on
$\hh^n$ is given by
$$
X_{2n+1}=\frac{\partial}{\partial t}, \hs
X_j=\frac{\partial}{\partial x_j}+2x_{n+j}\frac{\partial}{\partial
t}, \hs X_{n+j}= \frac{\partial}{\partial
x_{n+j}}-2x_j\frac{\partial}{\partial t},\ j=1,\,\cdots,\,n.$$ All
non-trivial commutators are $[X_j,\,X_{n+j}]=-4X_{2n+1}$,
$j=1,\,\cdots,\,n$. The sub-Laplacian has the form
$\Delta_{\hh^n}=\sum_{j=1}^{2n}X_j^2.$

Let $V$ be a nonnegative locally integrable function on $\hh^n$.
Define the sub-Laplacian Schr\"odinger operator by
$\cl\equiv-\Delta_{\hh^n}+V.$ Denote by $\{T_t\}_{t>0}\equiv
\{e^{-t\cl}\}_{t>0}$ the semigroup generated by $\cl$ and
by $\{\wz T_t\}_{t>0}\equiv
\{e^{t\Delta_{\hh^n}}\}_{t>0}$ the semigroup generated by $-\Delta_{\hh^n}$.

Let $V\in\cb_q({\hh^n}, d, dx)$ with $q\in(n+1, 2n+2]$ and $\rho$ be as in \eqref{2.4}.
Then $\{T_t\}_{t>0}$ and $\{Q_t\}_{t>0}$ satisfy Proposition \ref{p7.1}
with $d$ replaced by $2(n+2)$ and $|x-y|$ replaced by $d(x,\,y)$; see \cite{ll08}.

Observe that $\{\wz T_{t^2}\}_{t>0}$ is a continuous $(1,
N)$-$\ati$ for all positive constants $N$. Thus $\{T_{t^2}\}_{t>0}$ and
$\{\wz T_{t^2}\}_{t>0}$ satisfy the assumption \eqref{5.3}.
Moreover, the $L^2(\hh^n)$-boundedness of $g$-function was obtained in \cite{ll08}.
Using these facts and applying Theorems \ref{t5.1}, \ref{t5.2} and \ref{t6.1},
and Corollaries \ref{c5.1}, \ref{c5.2} and \ref{c6.1}, we have the following
conclusions.

\begin{prop}\hspace{-0.2cm}{\bf.}\hspace{0.2cm}\label{p7.4}
Let $q\in(n+1, \fz]$,  $V\in \cb_q(\hh^n,\,d,\,dx)$ and
$\rho $ be as in \eqref{2.4}.
Then there exists a positive constant $C$ such that
for all $f\in\bmo_{\rho}(\hh^n)$,
$T^+(f)$, $\wz T^+(f)$, $\wz T^+_\rho(f)$, $P^+(f)$, $\wz P^+(f)$,
 $\wz P^+_\rho(f),$ $g(f)$, $[g(f)]^2\in\blo_{\rho}(\hh^n)$ and
\begin{eqnarray*}
&&\|T^+(f)\|_{\blo_{\rho}(\hh^n)}
+\|\wz T^+(f)\|_{\blo_{\rho}(\hh^n)}+
\|T^+_\rho(f)\|_{\blo_{\rho}(\hh^n)}+\|P^+(f)\|_{\blo_{\rho}(\hh^n)}\\
&&\hs\hs+\|\wz P^+(f)\|_{\blo_{\rho}(\hh^n)}
+\|\wz P^+_\rho(f)\|_{\blo_{\rho}(\hh^n)}
+\|g(f)\|_{\blo_{\rho}(\hh^n)}+\|[g(f)]^2\|_{\blo_{\rho}(\hh^n)}^{1/2}\\
&&\hs\le C\|f\|_{\bmo_{\rho}(\hh^n)}.
\end{eqnarray*}
\end{prop}

We also point out that when $\rho$ is as in \eqref{2.4},
Lin and Liu \cite{ll08} introduced $\bmo_\rho(\hh^n)$ and
established the boundedness of
$T^+$, $P^+$ and $g$ on $\bmo_\rho(\hh^n)$.
The results in this subsection improve their corresponding results.

\medskip
\noindent 7.4. {\bf Schr\"odinger operators on connected
and simply connected nilpotent Lie groups.}
Let $\mathbb G$ be a connected and simply
connected nilpotent Lie group. Let
$X\equiv\{X_1,\cdots,X_k\}$ be left invariant vector fields on
$\bbg$ satisfying the H\"ormander condition that
$\{X_1,\,\cdots,\,X_k\}$ together with their commutators of order
$\le m$ generates the tangent space of $\bbg$ at each point of
$\bbg$. Let $d$ be the Carnot-Carath\'eodory (control) distance on
$\bbg$ associated to $\{X_1,\,\cdots,\,X_k\}$. Fix a left invariant
Haar measure $\mu$ on $\bbg$. Then for all $x\in \bbg$,
$V_r(x)=V_r(e)$; moreover, there exist $\kz$, $D\in(0, \fz)$ with
$\kz\le D$ such that for all $x\in\bbg$,
\begin{equation}\label{7.1}
C^{-1} r^\kz\le V_r(x)\le Cr^\kz
\end{equation}
when
$r\in(0, 1]$, and $C^{-1} r^D\le V_r(x)\le Cr^D$ when $r\in(1, \fz)$; see
\cite{nsw85} and \cite{v88}. Thus $({\mathbb
\bbg},\,d,\,\mu)$ is an RD-space.

The sub-Laplacian is given by $\Delta_\bbg\equiv\sum_{j=1}^kX_j^2.$
Denote by $\{\wz T_t\}_{t>0}\equiv\{e^{t\Delta_\bbg}\}_{t>0}$ the
semigroup generated by $-\Delta_\bbg$. Then there exist positive
constants $C,\,C_7$ and $\wz C_7$ such that for all $t\in(0, \fz)$ and
$x,\,y\in\bbg$,
\begin{equation}\label{7.2}
C^{-1}\frac1{V_{\sqrt t}(x)}\exp\lf\{-\frac{[d(x,\,y)]^2}{\wz
C_7t}\r\}\le \wz T_t(x,\,y)\le C\frac1{V_{\sqrt
t}(x)}\exp\lf\{-\frac{[d(x,\,y)]^2}{C_7t}\r\},
\end{equation}
that for all $t\in(0, \fz)$ and $x,\,y,\,y'\in\bbg$ with $d(y,\,y')\le
d(x,\,y)/4$,
\begin{equation}\label{7.3}
 |\wz T_t(x,\,y)-\wz T_t(x,\,y')|\le C\frac{d(y,\,y')}{\sqrt t}
\frac1{V_{\sqrt t}(x)}\exp\lf\{-\frac{[d(x,\,y)]^2}{C_7t}\r\},
\end{equation}
and that for all $t\in(0, \fz)$ and $x,\,y\in\bbg$,
\begin{equation*}
 \int_\bbg \wz T_t(x,\,z)\,d\mu(z)=1=\int_\bbg \wz T_t(z,\,y)\, d\mu(z);
\end{equation*}
see, for example, \cite{v88}.

Define the radial maximal operator $\wz T^+$ by $\wz
T^+(f)(x)\equiv\sup_{t>0}|\wz T_t(f)(x)|$ for all $x\in\bbg$. Then
by \eqref{7.2}, it is easy to see that
$\wz T^+$ is bounded on $L^p(\bbg)$ for $p\in(1,\,\fz]$.

Let $V$ be a nonnegative locally integrable function on $\bbg$. Then
the sub-Laplace Schr\"odinger operator $\cl$ is defined by
$\cl\equiv-\Delta_\bbg+V.$ The operator  $\cl$ generates a semigroup
of operators $\{T_t\}_{t>0}\equiv\{e^{-t\cl}\}_{t>0}$, whose kernels
are denoted by $\{T_t(x,\,y)\}_{t>0}$.
Define the radial maximal operator $T^+$ by
$T^+(f)(x)\equiv\sup_{t>0}|e^{-t\cl}(f)(x)|$ for all $x\in\bbg$. Then
from Lemma \ref{l7.1} below,
it is easy to see that $T^+$ is bounded on $L^p(\bbg)$ for $p\in(1,\,\fz]$.

Let $q>D/2$,  $V\in \cb_q(\bbg,\,d,\,\mu)$ and
$\rho$ be as in \eqref{2.4}. Then Li
\cite{l99} established some basic results concerning $\cl$, which
include estimates for fundamental solutions of $\cl$ and the
boundedness on Lebesgue spaces of some operators associated to
$\cl$. To apply the results obtained in Sections 5 and
6 to $\cl$, we need the following estimate, which
is a consequence of Proposition 5.2 and (5.12) in \cite{yz08}
together with the symmetry of $T_t$ and the fact that for
all $x,\,y\in\bbg$ and $t\in(0,\fz)$, $V_t(x)\sim V_t(y)$.
We omit the details.

\begin{lem}\hspace{-0.2cm}{\bf.}\hspace{0.2cm}\label{l7.1}
Let $q\in(D/2, D]$ and $V\in \cb_q(\bbg,\,d,\,\mu)$. Then for all $N\in(0, \fz)$,
there exist positive constants $C$ and $C_8$, where $C_8$
is independent of $N$, such that
for all $t\in(0, \fz)$ and $x,\,y\in\bbg$,
$$0\le T_t(x,\,y)\le C\frac1{V_{\sqrt t}(x)}
\exp\lf\{-\frac{[d(x,\,y)]^2}{C_8t}\r\}\lf[\frac{\rho(x)}{\rho(x)+\sqrt
t}\r]^N\lf[\frac{\rho(y)}{\rho(y)+\sqrt
t}\r]^N.$$
\end{lem}

For $t\in[0,\fz)$, set $E_t\equiv\wz T_t-T_t.$ Denote also by $E_t$ the
kernel of $E_t$.
The following estimate for $E_t$ was
established in \cite{yz08}.

\begin{lem}\hspace{-0.2cm}{\bf.}\hspace{0.2cm}\label{l7.2}
If $q\in(D/2, D]$  and $V\in \cb_q(\bbg,\,d,\,\mu)$, then
for all $N\in(0, \fz)$, there exist positive constants $C$ and $C_9$,
where $C_9$ is independent of $N$, such that
for all $t\in(0, \fz)$ and $x,\,y\in \bbg$,
\begin{equation*}
0\le E_t(x,\,y)\le
C\lf[\frac{\sqrt t}{\sqrt t+\rho(x)}\r]^{2-D/q} \frac1{V_{\sqrt t}(x)}
\exp\lf\{-\frac{[d(x,\,y)]^2}{C_9t}\r\}.
\end{equation*}
\end{lem}

Moreover, to estimate the regularity of $T_t$, we need the
regularity of $E_t$. To this end, we recall the following lemma.

\begin{lem}\hspace{-0.2cm}{\bf.}\hspace{0.2cm}\label{l7.3}
If $q\in(D/2, D]$ and $V\in \cb_q(\bbg,\,d,\,\mu)$, then for all positive constants
$C$ and $\wz C$, there exists positive constant $A_6$ such that for all
$x\in\bbg$ and $t>0$, when $\sqrt t<C\rho(x)$,
$$\int_\bbg\frac{V(z)}{V_{\sqrt t}(x)}
\exp\lf\{-\frac{[d(x,\,z)]^2}{\wz Ct}\r\}\,d\mu(z) \le A_6
\frac1t\lf[\frac{\sqrt t}{\rho(x)}\r]^{2-D/q},$$
while when $\sqrt t\ge C\rho(x)$,
$$\int_\bbg\frac{V(z)}{V_{\sqrt t}(x)}
\exp\lf\{-\frac{[d(x,\,z)]^2}{\wz Ct}\r\}\,d\mu(z) \le A_6
\frac1t\lf[\frac{\sqrt t}{\rho(x)}\r]^{\ell_0},$$
where $\ell_0$ is a positive constant independent of $C, \wz C$ and $A_6$.
\end{lem}

We remark that Lemma \ref{l7.3} with $\sqrt t<C\rho(x)$ is just
Lemma 5.1 of \cite{yz08}. For  $\sqrt t\ge C\rho(x)$, the result
can be proved similarly. We omit the details.

\begin{lem}\hspace{-0.2cm}{\bf.}\hspace{0.2cm}\label{l7.4}
If  $q\in(D/2, D]$  and $V\in \cb_q(\bbg,\,d,\,\mu)$, then for all
$\dz'\in(0, 2-D/q)$, there exist positive constants $C$ and
$A_7$, where $A_7$ is independent of $\dz$, such that for all
$t\in(0, \fz)$ and $x',\,x,\,y\in \bbg$ with
$d(x,x')<\min\{d(x,y)/4,\,\rho(x)\}$,
\begin{equation*}
|E_t(x,\,y)-E_t(x',\,y)|\le C\lf[\frac{d(x,x')}{\rho(y)}\r]^{\dz'}
\frac1{V_{\sqrt t}(x)}\exp\lf\{-\frac{[d(x,\,y)]^2}{A_7t}\r\}.
\end{equation*}
\end{lem}

\begin{pf}\rm
Let $x',\,x,\,y\in \bbg$ with
$d(x,x')<\min\{d(x,y)/4,\,\rho(x)\}$.
Notice that if $d(x,x')<d(x,y)/4$, then $d(x,y)\sim d(x',y)$.
We first prove that for all $\dz'\in(0, 2-D/q)$ and $x$, $y\in\bbg$,
\begin{equation}\label{7.4}
|E_t(x,\,y)-E_t(x',\,y)|\ls\lf[\frac{d(x,x')}{\rho(y)}\r]^{\dz'}
\frac1{V_{\sqrt t}(x)}.
\end{equation}

If $d(x,\,x')\ge\rho(y)$, then \eqref{7.4} follows
from Lemma \ref{l7.2}. If $d(x,\,x')<\rho(y)$ and $t\le
2[d(x,\,x')]^2$, another application of Lemma \ref{l7.2} together with
the symmetry of $T_t$ and $\wz T_t$ also yields \eqref{7.4}.
Thus we may assume that $d(x,\,x')<\rho(y)$ and $t> 2[d(x,\,x')]^2$.

Recall (see, for example, \cite{yz08,d05}) that for all $x$, $y\in\bbg$,
$$E_t(x,\,y)=\wz T_t(x,\,y)-T_t(x,\,y)= \int_0^t\int_\bbg \wz
T_{t-s}(x,\,z)V(z) T_s(z,\,y)\,d\mu(z)\,ds.$$
We write
\begin{eqnarray*}
|E_t(x,\,y)-E_t(x',\,y)|&&\le
\int_0^t\int_\bbg |\wz T_{t-s}(x,\,z)
-\wz T_{t-s}(x',\,z)|V(z)T_s(z,\,y)\,d\mu(z)\,ds\\
&&=\int_0^{t/2}\int_\bbg |\wz T_{t-s}(x,\,z)
-\wz T_{t-s}(x',\,z)|V(z) T_s(z,\,y)\,d\mu(z)\,ds\\
&&\hs+\int_0^{t/2}\int_\bbg |\wz T_s(x,\,z)-\wz T_s(x',\,z)
|V(z) T_{t-s}(z,\,y)\,d\mu(z)\,ds\\
&&\equiv {\rm F_1}+{\rm F_2}.
\end{eqnarray*}

To estimate ${\rm F_1}$, we consider the following two cases.
Case (i) $t<2[\rho(y)]^2$. For $s\in(0, t/2)$, we have $t-s\sim t$. By \eqref{7.3},
Lemma \ref{l7.1}, Lemma \ref{l7.3}, the symmetry of $\wz T_t$, the assumption that
$D/2<q\le D$ and the fact that
$V_{r}(x)\sim V_r(y)$ for all $x$, $y\in\bbg$ and $r\in(0,\fz)$,
we have
\begin{eqnarray*}
{\rm F_1}&&\ls \frac{d(x,x')}{\sqrt t}\frac1{V_{\sqrt t}(x)}
\int_0^{t/2}\int_\bbg V(z) \frac1{V_{\sqrt{s}}(y)}
\exp\lf\{-\frac{[d(z,\,y)]^2}{Cs}\r\}\,d\mu(z)\,ds\\
&&\ls \frac{d(x,x')}{\sqrt t}\frac1{V_{\sqrt t}(x)}
\int_0^{t/2}\frac1s\lf[\frac{\sqrt s}{\rho(y)}\r]^{2-D/q}\,ds
\ls \lf[\frac{d(x,x')}{\rho(y)}\r]^{2-D/q}\frac1{V_{\sqrt t}(x)}.
\end{eqnarray*}

 Case (ii) $t\ge 2[\rho(y)]^2$. Let $\ell_0$ be as in Lemma \ref{l7.3}
 and $N>\ell_0$.
 Using \eqref{7.3}, Lemma \ref{l7.1} and Lemma
\ref{l7.3},  we have
\begin{eqnarray*}
 {\rm F_1}&&\ls \frac{d(x,x')}{\sqrt t}\frac1{V_{\sqrt t}(x)}
\int_0^{t/2}\int_\bbg V(z) \frac1{V_{\sqrt{s}}(y)}
\exp\lf\{-\frac{[d(z,\,y)]^2}{Cs}\r\}
\lf[\frac{\rho(y)}{\sqrt{s}+\rho(y)}\r]^N\,d\mu(z)\,ds\\
&&\ls \frac{d(x,x')}{\sqrt t}\frac1{V_{\sqrt t}(x)}
\lf\{\int_0^{[\rho(y)]^2}\frac1s\lf[\frac{\sqrt
s}{\rho(y)}\r]^{2-D/q}\,ds+
\int_{[\rho(y)]^2}^{t/2}\frac1s\lf[\frac{\rho(y)}{\sqrt{s}}\r]^{N-\ell_0}\,ds\r\}\\
&&\ls\frac{d(x,x')}{\rho(x)}\frac1{V_{\sqrt t}(x)}.
\end{eqnarray*}

To estimate ${\rm F_2}$, we further write
\begin{eqnarray*}
{\rm F_2}&&\ls\int_0^{[d(x,x')]^2}\int_\bbg
|\wz T_{s}(x,\,z)-\wz T_{s}(x',\,z)|V(z)T_{t-s}(z,\,y)\,d\mu(z)\,ds\\
&&\hs+\int_{[d(x,x')]^2}^{t/2}\int_{W_1}\cdots
+\int_{[d(x,x')]^2}^{t/2}\int_{W_2}\cdots\equiv {\rm H_1}+{\rm H_2}+{\rm H_3},
\end{eqnarray*}
where $W_1\equiv\{z\in\bbg:\ d(x,x')>d(x,z)/4\}$ and $W_2\equiv\bbg\setminus
W_1$.

Since $d(x, x')\le\rho(x)$ together with \eqref{3.7}
implies that $\rho(x')\sim \rho(x)$, by
$d(x, z)\sim d(x',z)$, \eqref{7.2}, \eqref{7.3} and Lemma \ref{l7.3},
we obtain
\begin{eqnarray*}
{\rm H_1}&&\ls\frac1{V_{\sqrt{t}}(y)}\lf[\frac{\rho(y)}{\sqrt{t}+\rho(y)}\r]^N
\int_0^{[d(x,x')]^2}\int_\bbg \frac1{V_{\sqrt{s}}(x)}
\exp\lf\{-\frac{[d(z,\,x)]^2}{Cs}\r\}V(z) \,d\mu(z)\,ds\\
&&\ls\frac1{V_{\sqrt{t}}(y)}\lf[\frac{\rho(y)}{\sqrt{t}+\rho(y)}\r]^N
\lf[\frac{d(x,x')}{\rho(y)}\r]^{2-D/q}\lf[\frac{\rho(y)}{\rho(x)}\r]^{2-D/q}.
 \end{eqnarray*}

Recall (2.4) in \cite{yz08} that
for any fixed $y\in\bbg$ and $0<r<R<\fz$,
\begin{eqnarray*}
\frac{r^2}{V_r(y)}\dint_{B(y,\,r)}V(z)\,d\mu(z)
&&\ls\lf(\frac rR\r)^{2-D/q}\frac{R^2}{V_R(y)}\int_{B(y,\,R)}V(z)\,d\mu(z).
\end{eqnarray*}
Let $\kappa$ be as in \eqref{7.1}. By the inequality above,
the assumption that $q\in(D/2, D]$,
\eqref{7.3} and Lemma \ref{l7.1},
we have
 \begin{eqnarray*}
{\rm H_2}&&\ls\frac1{V_{\sqrt{t}}(y)}\lf[\frac{\rho(y)}{\sqrt{t}+\rho(y)}\r]^N
\int_{[d(x,x')]^2}^{t/2}\int_{W_1} \frac{d(x,x')}{\sqrt
s}\frac1{V_{\sqrt{s}}(x)}
V(z) \,d\mu(z)\,ds\\
&&\ls\frac1{V_{\sqrt{t}}(y)}\lf[\frac{\rho(y)}{\sqrt{t}+\rho(y)}\r]^N
\lf[\frac{d(x,x')}{\rho(y)}\r]^{2-D/q}\lf[\frac{\rho(y)}{\rho(x)}\r]^{2-D/q}.
 \end{eqnarray*}
For ${\rm H_3}$, if $t\le 2[\rho(x)]^2$, by $d(x,x')\le d(x,z)/2$, \eqref{7.3},
the assumption that $q\in(D/2, D]$, Lemma \ref{l7.1} and Lemma \ref{l7.3} with $\sqrt
s\le\sqrt2\rho(x)$, we obtain
 \begin{eqnarray*}
{\rm H_3}&&\ls\frac1{V_{\sqrt{t}}(y)}\lf[\frac{\rho(y)}{\sqrt{t}+\rho(y)}\r]^N
\int_{[d(x,x')]^2}^{t/2}\int_{W_2} \frac{d(x,x')}{\sqrt
s}\frac{V(z)}{V_{\sqrt{s}}(x)}\exp\lf\{-\frac{[d(z,\,x)]^2}{Cs}\r\}
\,d\mu(z)\,ds\\
&&\ls \frac1{V_{\sqrt{t}}(y)}\lf[\frac{\rho(y)}{\sqrt{t}+\rho(y)}\r]^N
\lf[\frac{d(x,x')}{\rho(y)}\r]^{\dz'}\lf[\frac{\rho(y)}{\rho(x)}\r]^{\dz'}.
 \end{eqnarray*}

If $t>2[\rho(x)]^2$, similarly to the above estimate, using Lemma
\ref{l7.3} with $\sqrt s>\rho(x)$, we have
\begin{eqnarray*}
{\rm H_3}&&\ls\frac1{V_{\sqrt{t}}(y)}\lf[\frac{\rho(y)}{\sqrt{t}+\rho(y)}\r]^N
\lf\{\int_{[d(x,x')]^2}^{[\rho(x)]^2}+\int_{[\rho(x)]^2}^{t/2}\r\}\\
&&\hs\times\int_{W_2} \frac{d(x,x')}{\sqrt
s}\frac1{V_{\sqrt{s}}(x)}\exp\lf\{-\frac{[d(z,\,x)]^2}{Cs}\r\}
V(z) \,d\mu(z)\,ds\\
&&\ls\frac1{V_{\sqrt{t}}(y)}\lf[\frac{\rho(y)}{\sqrt{t}+\rho(y)}\r]^{N-\ell_0+1}
\lf[\frac{d(x,x')}{\rho(y)}\r]^{2-D/q}
\lf[\frac{\rho(y)}{\rho(x)}\r]^{\ell_0}.
\end{eqnarray*}

Let $k_0$ be as in Definition \ref{d2.2}. Observing that
$$\frac{\rho(y)}{\rho(x)}\ls \lf[1+\frac{d(x,y)}{\sqrt t}
\frac{\sqrt t}{\rho(y)}\r]^{k_0}
 \ls \lf[1+\frac{\sqrt t}{\rho(y)}\r]^{k_0}\lf[1+\frac{d(x,y)}{\sqrt t}\r]^{k_0},$$
we obtain \eqref{7.4} by taking $N$ large enough. The estimate \eqref{7.4}
together with Lemma \ref{l7.2} implies the desired estimate, which
completes the proof of Lemma \ref{l7.4}.
\hfill$\Box$\end{pf}

\begin{lem}\hspace{-0.2cm}{\bf.}\hspace{0.2cm}\label{l7.5}
If $q\in(D/2, D]$ and $V\in \cb_q(\bbg,\,d,\,\mu)$,
then for all $N\in(0, \fz)$ and
$\dz'\in(0,\,2-D/q)$, there exist positive constants $\wz C$ and $C$,
where $C$ is independent of $N$, such that for all
$t\in(0, \fz)$ and $x',\,x,\,y\in \bbg$ with $d(x,x')<\sqrt t$,
\begin{eqnarray}\label{7.5}
|T_t(x,\,y)-T_t(x',\,y)|&&\le \wz C\lf[\frac{d(x,x')}{\sqrt t}\r]^{\dz'}
\frac1{V_{\sqrt t}(x)}\exp\lf\{-\frac{[d(x,\,y)]^2}{Ct}\r\}\nonumber\\
&&\hs\times\lf[\frac{\rho(y)}{\sqrt t+\rho(y)}\r]^{N}\lf[\frac{\rho(x)}{\sqrt
t+\rho(x)}\r]^{N}.
\end{eqnarray}
\end{lem}

\begin{pf}\rm
Let $x',\,x,\,y\in \bbg$ with $d(x,x')<\sqrt t$.
We first consider the case $d(x,x')\le d(x,y)/4$. Using \eqref{2.3}, we obtain
\begin{equation}\label{7.6}
\frac{\rho(x')}{\sqrt t+\rho(x')}\ls \lf[\frac{\rho(x)}{\sqrt
t+\rho(x)}\r]^{1/(1+k_0)} \lf[1+\frac{d(x,y)}{\sqrt
t}\r]^{k_0/(1+k_0)}.
\end{equation}
Since $d(x,y)\sim d(x',y)$, by Lemma \ref{l7.1}, we have
\begin{equation}\label{7.7}
|T_t(x,\,y)-T_t(x',\,y)|\ls \frac1{V_{\sqrt
t}(x)}\exp\lf\{-\frac{[d(x,\,y)]^2}{Ct}\r\} \lf[\frac{\rho(y)}{\sqrt
t+\rho(y)}\r]^{N}\lf[\frac{\rho(x)}{\sqrt t+\rho(x)}\r]^{N}.
\end{equation}
If $d(x,x')\ge \rho(x)$, then \eqref{7.5} follows from
\eqref{7.7}. If $d(x,x')< \rho(x)$, then by Lemma \ref{l7.4}, we obtain
$$|E_t(x,\,y)-E_t(x',\,y)|\ls\lf[\frac{d(x,x')}{\sqrt t}\r]^{\dz'}
\lf[\frac{\sqrt t}{\rho(y)}\r]^{\dz'} \frac1{V_{\sqrt
t}(x)}\exp\lf\{-\frac{[d(x,\,y)]^2}{Ct}\r\}.$$
This together with \eqref{7.3} and $T_t=\wz
T_t-E_t$ gives that
$$
|T_t(x,\,y)-T_t(x',\,y)|\ls\lf[\frac{d(x,x')}{\sqrt t}\r]^{\dz'}
\lf[1+\frac{\sqrt t}{\rho(y)}\r]^{\dz'} \frac1{V_{\sqrt
t}(x)}\exp\lf\{-\frac{[d(x,\,y)]^2}{Ct}\r\}.
$$
Then \eqref{7.5} follows from this and \eqref{7.7}.

Now we assume that $d(x,x')\ge d(x,y)/4$. In this case, $d(x,y)<4\sqrt{t}$.
Write
\begin{eqnarray*}
 |T_t(x,\,y)-T_t(x',\,y)|&&\ls\int_\bbg  |T_{t/2}(x,\,z)-T_{t/2}(x',\,z)
 |T_{t/2}(z,\,y)\,d\mu(z)\\
&&\ls\int_{W_1} |T_{t/2}(x,\,z)-T_{t/2}(x',\,z)|T_{t/2}(z,\,y)\,d\mu(z)
+\int_{W_2}\cdots\equiv {\rm I_1}+{\rm I_2},
\end{eqnarray*}
where $W_1$ and $W_2$ are as in Lemma \ref{l7.4}.

By Lemma \ref{l7.1}, we have
\begin{eqnarray*}
{\rm I_1}&&\ls\lf[\frac{V(x,x')}{V_{\sqrt t}(x)}+\frac{V(x,x')}{V_{\sqrt t}(x')}\r]
\frac1{V_{\sqrt t}(y)}\lf[\frac{\rho(y)}{\sqrt t+\rho(y)}\r]^{N}
\ls\lf[\frac{d(x,x')}{\sqrt t}\r]^{\dz''}\frac1{V_{\sqrt
t}(y)}\lf[\frac{\rho(y)}{\sqrt t+\rho(y)}\r]^{N},
\end{eqnarray*}
where $\dz''\in(\dz',2-D/q)$.

Using \eqref{7.5} with $d(x,x')\le d(x,z)/4$ and Lemma \ref{l7.1}, we obtain
\begin{eqnarray*}
{\rm I_2}\ls \lf[\frac{d(x,x')}{\sqrt t}\r]^{\dz''}
\frac1{V_{\sqrt t}(x)}\int_{W_2}T_{t/2}(z,\,y)\,d\mu(z)
\ls\lf[\frac{d(x,x')}{\sqrt t}\r]^{\dz''}\frac1{V_{\sqrt
t}(x)}\lf[\frac{\rho(y)}{\sqrt t+\rho(y)}\r]^{N}.
\end{eqnarray*}
By $d(x,y)\le 4 \sqrt t$, \eqref{7.6}
and Lemma \ref{l7.1}, we have \eqref{7.5} for $d(x,x')\ge d(x,y)/4$,
which completes the proof of Lemma \ref{l7.5}.
\hfill$\Box$\end{pf}

For all $x,\,y\in\bbg$ and $t\in(0, \fz)$, define
$$Q_t(x,\,y)\equiv t^2\frac{d}{ds}\Big|_{s=t^2}T_s(x,y).$$
Following the proof of Proposition 4 in
\cite{dgmtz05}, we have the following result. We omit the details.

\begin{lem}\hspace{-0.2cm}{\bf.}\hspace{0.2cm}\label{l7.6}
Let $q\in (D/2, D]$, $\bz\in(0,\,2-D/q)$ and $N\in\nn$.
There exist positive constants $\wz C$ and $C$,
where  $C$ is independent of $N$, such that
for all $t\in(0,\,\fz)$ and $x,\,x',\,y\in\bbg$ with $d(x,\,x')\le
\frac t2$,

(i) $|Q_t(x,\,y)|\le \wz C\frac1{V_t(x)}\exp \{-\frac
{[d(x,\,y)]^2}{Ct^2}\}[\frac {\rho(x)}{t+\rho(x)}]^N
[\frac {\rho(y)}{t+\rho(y)}]^N$;

(ii)  $|Q_t(x,\,y)-Q_t(x',\,y)|\le\wz C
 [\frac {d(x,\,x')}{t}]^{\bz}\frac1{V_t(x)}\exp \{-\frac
{[d(x,\,y)]^2}{Ct^2}\}[\frac {\rho(x)}{t+\rho(x)}]^N
[\frac {\rho(y)}{t+\rho(y)}]^N$;

(iii) $|\int_\cx Q_t(x,\,y)d\mu(y)|\le\wz C
[\frac t{\rho(x)}]^{2-D/q} [\frac{\rho(x)}{t+\rho(x)}]^N$.
\end{lem}

\begin{rem}\hspace{-0.2cm}{\bf.}\hspace{0.2cm}\label{r7.1}\rm
Let $q_1$, $q_2\in(D/2,\fz]$ with $q_1<q_2$.
Recall that $\cb_{q_2}(\bbg)\subset \cb_{q_1}(\bbg)$.
Therefore, Lemmas \ref{l7.1}  through 7.6 hold for all
$q\in(D/2, \fz]$.
\end{rem}

Observe that $\{\wz T_{t^2}\}_{t>0}$ is a continuous $(1,
N)$-$\ati$ for all positive constants $N$. Thus $\{T_{t^2}\}_{t>0}$ and
$\{\wz T_{t^2}\}_{t>0}$ satisfy the assumption \eqref{5.3}.
Moreover, the $L^2(\bbg)$-boundedness of $g$-function
can be obtained by the same argument
as in Lemma 3 of \cite{dgmtz05}.
Using these facts and applying Theorems \ref{t5.1}, \ref{t5.2} and \ref{t6.1},
and Corollaries \ref{c5.1}, \ref{c5.2} and \ref{c6.1},
we have the following conclusions.

\begin{prop}\hspace{-0.2cm}{\bf.}\hspace{0.2cm}\label{p7.5}
Let $q\in(D/2, \fz]$,  $V\in \cb_q(\bbg,\,d,\,\mu)$ and
$\rho$ be as in \eqref{2.4}.
There exists a positive constant $C$ such that
for all $f\in\bmo_{\rho}(\bbg)$,
$T^+(f)$, $\wz T^+(f)$, $\wz T^+_\rho(f)$, $P^+(f)$, $\wz P^+(f)$,
 $\wz P^+_\rho(f),$ $g(f),\,[g(f)]^2\in\blo_{\rho}(\bbg)$ and
\begin{eqnarray*}
&&\|T^+(f)\|_{\blo_{\rho}(\bbg)}
+\|\wz T^+(f)\|_{\blo_{\rho}(\bbg)}+
\|T^+_\rho(f)\|_{\blo_{\rho}(\bbg)}+\|P^+(f)\|_{\blo_{\rho}(\bbg)}\\
&&\hs\hs+\|\wz P^+(f)\|_{\blo_{\rho}(\bbg)}
+\|\wz P^+_\rho(f)\|_{\blo_{\rho}(\bbg)}
+\|g(f)\|_{\blo_{\rho}(\bbg)}
+\|[g(f)]^2\|_{\blo_{\rho}(\bbg)}^{1/2}\\
&&\hs\le C\|f\|_{\bmo_{\rho}(\bbg)}.
\end{eqnarray*}
\end{prop}

\bigskip

\noindent{\bf Acknowledgements.}  The authors would like to thank the referees for
their several valuable remarks which improve the presentation of this
article.

\begin{center}
\section*{REFERENCES}
\end{center}

\medskip

\noindent{\it E-mail address}: {dcyang@bnu.edu.cn}

\noindent{\it E-mail address}: {dyyang@mail.bnu.edu.cn}

\noindent{\it E-mail address}: {yuanzhou@mail.bnu.edu.cn}
\end{document}